\documentclass[smallextended,numbook,runningheads,12pt]
{svjour3}
\smartqed
\usepackage{mathptmx}
\usepackage{amssymb}
\usepackage{mathrsfs}
\usepackage{amsfonts}
\usepackage{graphicx}
\usepackage{color}
\usepackage{amsmath}
\usepackage{psfrag}
\usepackage{booktabs}
\numberwithin{equation}{section}
\usepackage[top=1.2in, bottom=1.2in, left=1.3in,
right=1.3in]{geometry}
\usepackage{listings}
\numberwithin{equation}{section}

\newcommand{\R}{\mathbb{R}}
\newcommand{\N}{\mathbb{N}}
\newcommand{\E}{\mathbb{E}}

\newcommand{\BF}{\mathbf{F}}
\newcommand{\BG}{\mathbf{G}}

\newtheorem{thm}{Theorem}[section]

\newtheorem{prop}{Proposition}[section]
\newtheorem{cor}[thm]{Corollary}

\newtheorem{assumption}[thm]{Assumption}

\journalname{Journal of Scientific Computing}
\begin{document}

\title{An exponential integrator scheme for
time discretization of nonlinear stochastic wave equation
}
       \author{
        Xiaojie Wang }

\institute{Xiaojie Wang (Corresponding author) \at School of Mathematics and Statistics and School of Geosciences and Info-Physics,
       Central South University, {\small Changsha 410083, Hunan,  PR China }\\
       \email{x.j.wang7@csu.edu.cn, x.j.wang7@gmail.com}}

\date{Received: date / Accepted: date}

\maketitle

\begin{abstract}

This work is devoted to convergence analysis of an
exponential integrator scheme for semi-discretization
in time of nonlinear stochastic wave equation.
A unified framework is first set forth, which covers
important cases of additive and multiplicative noises.
Within this framework, the proposed scheme is shown to
converge uniformly in the strong $L^p$-sense with
precise convergence rates given. The abstract
results are then applied to several concrete examples.
Further, weak convergence rates of the scheme are examined for the case of additive noise. To analyze the weak error for the nonlinear case, techniques based on the Malliavin calculus were usually exploited in the literature. Under certain appropriate assumptions on the nonlinearity, this paper provides a weak error analysis,  which does not rely on the Malliavin calculus. The rates of weak convergence can, as expected, be improved in comparison with the strong rates.  Both strong and weak convergence results obtained here show that the proposed method achieves higher convergence rates than the implicit Euler and Crank-Nicolson time discretizations. Numerical results are finally reported to confirm our theoretical findings.

\keywords{nonlinear stochastic wave equation \and multiplicative noise \and exponential Euler scheme \and strong convergence \and weak convergence}

\subclass{60H35 \and 60H15 \and 65C30}

\end{abstract}

\section{Introduction}
\label{sect:intro}

In the present paper, we consider stochastic evolution equation of It\^{o} type in a separable Hilbert space $(U,\, \langle \cdot, \cdot \rangle_U, \,\|\cdot \|_U)$, given by
\begin{align}\label{eq:Abstr.SWE}
\left\{\!
    \begin{array}{ll}
    \mbox{d} \dot{u}(t) = -\Lambda u(t) \, \mbox{d}t + F(u(t))\, \mbox{d}t + G(u(t))\,\mbox{d}W(t), \quad t \in (0, T],
    \\
     u(0) = u_0, \, \dot{u}(0) = v_0,
    \end{array}\right.
\end{align}
where $T \in (0, \infty)$ and $\Lambda \colon \mathcal{D}(\Lambda)
\subset U \rightarrow U $ is a densely defined, linear unbounded,
positive self-adjoint operator with compact inverse (e.g., $\Lambda = - \Delta$ with homogeneous Dirichlet boundary condition). Under this assumption, there exists an increasing sequence of real numbers $\{\lambda_i\}_{i =1}^{\infty}$ and an orthonormal basis $\{e_i\}_{i =1}^{\infty}$ such that $\Lambda e_i = \lambda_i e_i$ and
\begin{equation}
0 < \lambda_1 \leq \lambda_2 \leq \cdots \leq \lambda_n (\rightarrow \infty).
\end{equation}
In Equation \eqref{eq:Abstr.SWE}, $\{u(t)\}_{t\in [0, T]}$ is regarded as a $U$-valued stochastic process and $\{\dot{u}(t)\}_{t\in [0, T]}$ stands for the time derivative of $u$. Moreover, $F \colon U \rightarrow U$, $G(u) \colon U \rightarrow U$ for $u \in U$ are deterministic mappings and $\{W(t)\}_{t \in[0, T]}$ is a cylindrical $Q$-Wiener process on a given probability space $\left(\Omega,\mathcal {F},\mathbb{P}\right)$ with normal filtration $\{\mathcal {F}_t\}_{t\in [0, T]}$.

The abstract equation \eqref{eq:Abstr.SWE} includes many stochastic wave equations (SWEs) in applications \cite{CE72,DKMNX09,Walsh86}. Since their true solutions  are rarely known explicitly, numerical simulations are often used to understand the behavior of the solutions. To do this, one has to discretize both the time interval $[0, T]$ and the infinite dimensional space $U$. As pointed out in \cite{PJ01}, the main difficulty in studying numerical schemes for stochastic partial differential equations (SPDEs) of evolutionary type lies in the treatment of the time discretization. On the one hand, development of effective high order time-stepping schemes are important since better time integration technology gives significant performance improvements in the numerical solution of SPDEs. This issue is, however, essentially difficult, even in the context of finite dimensional stochastic ordinary differential equations (SODEs) \cite{KP92}. Time discretizations of SPDEs encounter all the difficulties that arise in the time approximations of  both deterministic PDEs and finite dimensional SODEs as well as many more due to the infinite dimensional nature of the driving noise processes. On the other hand, as one can see in our forthcoming work, the arguments used in the following analysis of pure time discretization can be extended to the analysis of fully discretized scheme when combined with arguments used in the deterministic theory. For instance, optimal error estimates for the deterministic wave equation obtained in \cite{KLS10} allow us to analyze a fully discretized scheme with finite element spatial discretization. Therefore, and also for the sake of simplicity,  we first concentrate on the time discretization of \eqref{eq:Abstr.SWE} in this paper.


On the interval $[0, T]$, we construct a
uniform mesh $\mathcal{T}_M = \{ t_0, t_1, \cdots, t_M \}$ satisfying
\begin{equation}
  \mathcal{T}_M : {0 = t_0 < t_1 < t_2 < \cdots < t_M=T}
\end{equation}
with $\tau = T/M$, $M \in \N$,  being the time stepsize. This article is concerned with the following scheme for the time discretization of \eqref{eq:Abstr.SWE}:
\begin{align}\label{eq:exp.Euler}
\left\{\!
    \begin{array}{l}
    u_{m+1}= C(\tau) u_m + \Lambda^{-\frac{1}{2}} S(\tau)v_m + \tau \Lambda^{-\frac{1}{2}} S(\tau) F(u_m) + \Lambda^{-\frac{1}{2}} S(\tau) G(u_m)  \Delta W_m,
    \\
    v_{m+1}= -\Lambda^{\frac{1}{2}} S(\tau) u_m + C(\tau)v_m +  \tau C(\tau) F(u_m)
    +  C(\tau) G(u_m)  \Delta W_m,
    \end{array}\right.
\end{align}
where $C(t)\!=\! \cos(t\Lambda^{\frac{1}{2}})$
and $S(t)\!=\! \sin(t\Lambda^{\frac{1}{2}})$ for
$t \in [0, T]$ are the cosine and sine operators
and $\Delta W_m = W(t_{m+1}) - W(t_m)$ is the
Wiener increment. Here $u_{m}$ and $v_{m}$ are,
respectively, temporal approximations of $u(t)$
and $\dot{u}(t)$ at the grid points
$t_m \in\mathcal{ T}_M$. Rewriting \eqref{eq:exp.Euler}
in another abstract form (see \eqref{eq:Abstr.Exp.Euler}),
one can observe that the scheme \eqref{eq:exp.Euler}
can be identified with a version of the stochastic
exponential Euler scheme in the literature
\cite{JK09a,JKW11,KLNS11,LR04,LT10},
where the schemes were designed for parabolic SPDEs.
By choosing particular filter functions,
the stochastic exponential Euler scheme can be
viewed as a stochastic trigonometric method
for nonlinear problems \cite[Example 2]{CS12}.
Recently, Cohen, Larsson and Sigg \cite{CLS13}
applied the stochastic trigonometric method to
the linear SWE with additive noise, which coincides with
the scheme \eqref{eq:exp.Euler} when applied to the
linear case.
In the first part of this article, we establish the uniform $L^p$-convergence of the scheme \eqref{eq:exp.Euler} in a general framework for nonlinear problems, which of course extends the $L^2$-convergence results in \cite{CLS13} for the linear case.

To get started, we make the following assumptions.
\begin{assumption}\label{assump:FG}
Assume $F \colon U \rightarrow U $ and $G(u) \colon U \rightarrow U$ for $u \in U$ are deterministic mappings satisfying
\begin{align}\label{eq:Nf.condition1}
  \| F(u) \|_U + \big\| \Lambda^{\frac{\delta -1 }{2}} G(u) \big\|_{\mathcal{L}_2^0} \leq& L(\|u\|_U+1),
  \\
  \label{eq:Nf.condition2}
  \| F(u_1)-F(u_2) \|_U + \big\| \Lambda^{\frac{ \delta - 1}{2}} \big( G(u_1)-G(u_2) \big) \big\|_{\mathcal{L}_2^0} \leq& L \|u_1-u_2\|_U
\end{align}
for all $u,u_1, u_2 \in U$ and some positive constants $\delta, L \in (0, \infty)$.
\end{assumption}
For notations, we refer to the next section.
Throughout this article, we choose $L \geq 1$ to be a
universal constant in all assumptions.
Under Assumption \ref{assump:FG}, \eqref{eq:Abstr.SWE} possesses a unique mild solution (Theorem \ref{thm:unique.mild}). Moreover, the spatial regularity properties of the mild solution and its numerical approximations are investigated under the above assumption (Proposition \ref{prop:MB.X.Xm}). Further, we establish the main strong convergence result, Theorem \ref{thm:strong.conv.result}, which shows that the scheme \eqref{eq:exp.Euler} is convergent uniformly in the strong $L^p$-sense, $p \in [2, \infty)$. In particular, we show
\begin{equation}\label{eq:introd.mainstrongconv}
\Big\| \sup_{t_m \in \mathcal{T}_M } \big\| u(t_m) - u_m \big\|_U \Big\|_{L^p(\Omega, \R)} \leq C \tau^{\min(\delta, 1)}.
\end{equation}
Here, and throughout this work, $C$ is a constant that may vary from one place to another and depends on $T, \delta, L, p, \lambda_1$ and the initial data $u_0,v_0$, but is independent of $M$. To derive the uniform bound \eqref{eq:introd.mainstrongconv}, a continuous-time extension of $u_m$ is introduced and the group property of the operator family $E(t) = e^{t A}$ is exploited.
Then, several concrete examples including SWEs with  general additive noise, multiplicative trace class noise and multiplicative space-time white noise are presented to fit in the above setting. For the three different cases, \eqref{eq:introd.mainstrongconv} is then applied with precise convergence rates given. It is shown that the spatial regularity of the noise term determines the order of strong convergence. Moreover, unlike the parabolic case, the exponential Euler scheme achieves higher strong order than the implicit Euler and Crank-Nicolson time discretizations do (see also \cite{KLL12}). Another interesting finding is that,  for the particular case of multiplicative trace-class noise the scheme \eqref{eq:exp.Euler} can achieve the strong convergence order of one, which is, as one can see later, the ultimate limit on the rate one can achieve for \eqref{eq:exp.Euler}. This upper limit on convergence order can be also found in \cite{KLL12}.



For the case of additive noise satisfying
\begin{equation}\label{eq:introd.weak.QA}
\big\| \Lambda^{ -\frac{1}{2} + \beta} Q  \Lambda^{ -\frac{1}{2} } \big\|_{\mathcal{L}_1(U)} < \infty,
\end{equation}
we also measure the weak error $\big|\,\E \!\left[ \varphi (u(T))\right] - \mathbb{E} \! \left[\varphi (u_M)\right] \big|$ for a class of test functions $\varphi \colon U \rightarrow \R$ with $\varphi \in \mathcal{C}_b^2(U, \R)$. We mention that the essential condition \eqref{eq:introd.weak.QA} has been used in \cite{KLL11,KLL12} to carry out the weak convergence analysis of numerical schemes for the linear SWE. Also, we emphasize that the presence of the nonlinear term leads to non-trivial technical difficulties in the analysis of weak error. Indeed, in the weak error analysis for linear SPDEs with additive noise \cite{DP09,KLL11,KLL12,QY13}, whose solution can be written down explicitly, one can get rid of the irregular term involving the unbounded operator by a transformation of variables. This transformation, however, does not work for the nonlinear heat equation. Since the operator family $E(t) = e^{t A}$ is a group (see the next section), one can adapt the basic line of \cite{DD06} to make the weak error analysis easier. Instead of imposing very strong spatial regularity conditions on the nonlinearities \cite[Hypothesis 5.7, 5.8]{DD06}, we make certain relatively mild assumptions on the nonlinearity $F$ (c.f. Theorem \ref{thm:weak.conv.result} and Corollary \ref{cor:weak.error}). For instance, the assumptions do not require $F$ to be twice Fr\'{e}chet differentiable in $U$, which is in general not fulfilled for Nemytskij operators but was commonly demanded in the weak error analysis \cite{AL13,DD06,DA10,WG13a}. Moreover, a key condition on $F$ (see \eqref{eq:assumption.F'.weak})  is used to avoid techniques involving the Malliavin calculus. Under such assumptions we obtain
\begin{equation}
\big|\mathbb{E} \! \left[ \varphi (u(T))\right] - \mathbb{E} \! \left[ \varphi (u_M) \right] \big| \leq
C_{\epsilon} \tau^{\min( 2\beta, \frac{1}{2} + \beta -
\epsilon, 1) }
\end{equation}
with arbitrarily small $\epsilon >0$.
As expected, the rate of weak convergence is, up to an arbitrarily small $\epsilon >0$ in some cases, twice that of strong convergence. Eventually, the assumptions we put are illustrated by a concrete setting.


In summary, the contribution of this article to the numerical analysis of SWEs is twofold. On the one hand, a unified abstract framework is formulated and a high order time-stepping scheme is proposed for nonlinear SWEs. Within the abstract framework, strong convergence of the exponential Euler scheme is studied thoroughly. The framework is general as it covers cases of additive noise, multiplicative trace-class noise as well as multiplicative space-time white noise. On the other hand, we provide a weak error analysis which does not rely on the Malliavin calculus, while techniques based on the Malliavin calculus were usually exploited in the course of weak error analysis for the nonlinear case \cite{AL13,AKL13,Brehier12,BrehierKopec13,DA10,WG13a}. Both strong and weak convergence results show that the proposed scheme attains higher convergence rates than the implicit Euler and Crank-Nicolson time discretizations.


Finally,  we mention that parabolic SPDEs have been extensively studied from the numerical point of view. The readers are referred to a review article
\cite{JK09}, a relevant monograph \cite{JK11} and more recent doctor thesis \cite{Kruse12,Lindgren12}
for numerous references.
The numerical research of stochastic wave equation,
however, is in its beginning. Available results
\cite{CL07,CLS13,HE10,KLS10,KLL11,KLL12,MPW03,QS06,Walsh06,WGT13}
are much less compared with the numerical analysis of stochastic parabolic problems, which also
partly motivates this work.

A brief outline of this paper is as follows.
In the next section, we collect some preliminaries and
formulate an abstract framework.
In Section \ref{sect:stongconvergence}, we analyze the strong approximation error arising from the time discretization. Then  in Section \ref{sec:examples}, several examples are included,
which fit in the abstract setting, to illustrate our abstract results. Weak convergence of the scheme is studied in Section \ref{sec:weak.convergence} for the case of additive noise. Numerical results are presented at the end of this article.


\section{Preliminaries and framework}
\label{sect:framework}

Let $(U,\, \langle \cdot, \cdot \rangle_U, \,\|\cdot \|_U)$ and $(H,\, \langle \cdot, \cdot \rangle_H, \, \|\cdot \|_H)$
be two separable Hilbert spaces.
By $\mathcal{L}(U,H)$ we denote the space of bounded linear operators from $U$ to $H$ with the usual operator norm $\| \cdot \|_{\mathcal{L}(U, H)}$ and write $\mathcal{L}(U)=\mathcal{L}(U,U)$ to lighten the notation. Additionally, we need  spaces of nuclear and Hilbert-Schmidt operators \cite{DZ92,PR07}. The space of nuclear operators from $U$ to $H$ is denoted by $\mathcal{L}_1(U,H)$ and we write $\mathcal{L}_1(U) = \mathcal{L}_1(U,U)$.
If $\Gamma \in \mathcal{L}_1(U)$ is nonnegative and symmetric, then
\begin{align}
\|\Gamma\|_{\mathcal{L}_1(U)} = \text{Tr}(\Gamma):= \sum_{i=1}^{\infty} \langle \Gamma \psi_i, \psi_i \rangle_U,
\end{align}
where $\{\psi_i\}_{i \in \mathbb{N}}$ is an orthonormal basis of $U$ and the trace of a nuclear operator, namely, $\mbox{Tr}(\Gamma)$ for $\Gamma \in \mathcal{L}_1(U)$, is independent of the particular choice of the basis $\{\psi_i\}_{i\in \mathbb{N}}$. By $\N := \{1,2,3, \ldots\}$ we denote the natural numbers
and by $\mathcal{L}_2(U,H)$ we denote the space of Hilbert-Schmidt operators from $U$ to $H$, equipped with the norm
\begin{equation}\label{eq:normdef} \|\Gamma\|_{\mathcal{L}_2(U,H)} = \Big(\sum_{i=1}^{\infty}\|\Gamma \psi_i\|_H^2\Big)^{1/2},
\end{equation}
also not depending on the particular choice of the basis. Analogously, we write $\mathcal{L}_2(U) = \mathcal{L}_2(U,U)$. If $\Gamma_1\in \mathcal{L}(U,H)$ and $\Gamma_2\in \mathcal{L}_j(U)$, $j =1,2$, then $\Gamma_1 \Gamma_2 \in \mathcal{L}_j(U,H)$ for $j =1,2$, and
\begin{equation}\label{eq:LHS}
\| \Gamma_1 \Gamma_2\|_{\mathcal{L}_j(U,H)} \leq \|\Gamma_1\|_{\mathcal{L}(U,H)} \cdot \|\Gamma_2\|_{\mathcal{L}_j(U)}, \quad  j=1,2.
\end{equation}
Moreover, if $\Gamma_1 \in \mathcal{L}(H)$ and $\Gamma_2 \in \mathcal{L}_1(H)$, then both $\Gamma_1\Gamma_2$ and $\Gamma_2\Gamma_1$ belong to $\mathcal{L}_1(H)$ and
\begin{align} \label{ineq4:L1L2}
\text{Tr} (\Gamma_1\Gamma_2) = \text{Tr}(\Gamma_2\Gamma_1).
\end{align}

Let $\left(\Omega,\mathcal {F},\mathbb{P}\right)$ be a probability space with a normal filtration $\{\mathcal{F}_t\}_{0\leq t\leq T}$ and by $L^p(\Omega, U)$ we denote the space of $U$-valued integrable random variables with the norm defined by
$
\|X\|_{L^p(\Omega, U)} = \big(\mathbb{E} \!\left[\|X\|_U^p \right] \big)^{\frac{1}{p}}
$
for $p\in [2, \infty)$.
Furthermore, the driven stochastic process $W(t)$ in \eqref{eq:Abstr.SWE} is assumed to be a cylindrical $Q$-Wiener process in the stochastic basis $\left(\Omega,\mathcal {F},\mathbb{P}, \{\mathcal {F}_t\}_{0\leq t\leq T}\right)$,
with a covariance operator $Q:\! U \rightarrow U$,
which can be represented as follows \cite{DZ92,PR07}:
\begin{align}\label{W.representation}
W(t) = \sum_{i=1}^{\infty} \sqrt{q_i} \, \beta_i(t) \phi_i, \quad t \in [0, T],
\end{align}
where $\{\beta_i(t) \}_{i\in \{ n \in \N \colon q_n >0 \} }$ are a family of mutually independent real Brownian motions and $\{\phi_i\}_{i \in \mathbb{N}}$ form an orthonormal basis of $U$ consisting of eigenfunctions of $Q$ with $Q \phi_i = q_i \phi_i, \, q_i \geq 0, \, i \in \mathbb{N}$. The covariance operator $Q\in \mathcal{L}(U)$ is nonnegative and symmetric, but not necessarily of finite trace. As a result, the series in \eqref{W.representation} may not converge in $U$, but in some space $U_1$ into which $U$ can be embedded, see, e.g., \cite{DZ92,PR07}. Let $Q^{\frac{1}{2}}$ denote the unique positive square root of $Q$. Now we are able to introduce the separable Hilbert space $U_0 := Q^{\frac{1}{2}}(U)$ endowed with the inner product $\langle u, \hat{u} \rangle_0 = \langle Q^{-\frac{1}{2}} u, Q^{-\frac{1}{2}} \hat{u} \rangle_U$ for $u, \hat{u} \in U$, where $Q^{-\frac{1}{2}}$ is the pseudo inverse of $Q^{\frac{1}{2}}$ in the case when $Q$ is not one-to-one. For lighter notation, we use $\mathcal{L}_2^0$ to denote the separable Hilbert spaces $\mathcal{L}_2(U_0, U)$ and also $\mathcal{L}_2(U_0, H)$ when it causes no confusion.

At this point, we are ready to discuss the existence and uniqueness of the so-called mild solution of \eqref{eq:Abstr.SWE}. To this end, we rewrite  \eqref{eq:Abstr.SWE} as a stochastic evolution equation of first order in a product space $H$. Specifically, introducing a new variable $v = \dot{u}$ transforms \eqref{eq:Abstr.SWE} into the following Cauchy problem
\begin{equation}\label{eq:SEE}
\left\{
    \begin{array}{ll}
    \mbox{d} X(t) = A X(t) \mbox{d}t + \BF(X(t))\, \mbox{d}t + \BG(X(t))\,\mbox{d}W(t), \quad t \in (0, T],
    \\
     X(0)= X_0,
    \end{array}\right.
\end{equation}
where we denote
\begin{equation}\label{eq:SEE.notations}
X_0 = \left[ \!\begin{array}{c}
    u_0
    \\
    v_0
    \end{array}\!\right], \:
X = \left[ \!\begin{array}{c}
    u
    \\
    v
    \end{array}\!\right], \:
    A = \left[\! \begin{array}{cc}
    0 & I
    \\
    -\Lambda & 0
    \end{array}\!\right],\:
    \BF(X) = \left[\! \begin{array}{c}
    0
    \\
    F(u)
    \end{array}\!\right],\:
    \BG(X) = \left[\! \begin{array}{c}
    0
    \\
    G(u)
    \end{array}\!\right].
\end{equation}
Here and below by $I$ we mean the identity operator in $U$. In this way, we transfer the existence and uniqueness of the mild solution of \eqref{eq:Abstr.SWE} to the same problem for \eqref{eq:SEE}. Before proceeding further,
we need additional spaces and notations.
We introduce the Hilbert space $\dot{H}^{\gamma} = \mathcal{D}(\Lambda^{\frac{\gamma}{2}})$ for $\gamma \in \R$ \cite{Kruse12}, equipped with the inner product
\begin{equation}\label{eq:Lambda.norm}
  \langle u, \omega \rangle_{\dot{H}^{\gamma}} := \langle \Lambda^{\frac{\gamma}{2}} u, \Lambda^{\frac{\gamma}{2}} \omega \rangle_U = \sum_{i=1}^{\infty} \lambda_i^{\gamma} \langle u, e_i\rangle_U \langle \omega, e_i\rangle_U, \quad \gamma \in \mathbb{R},
\end{equation}
and the corresponding norm $\|u\|^2_{\gamma} = \langle u, u \rangle _{\dot{H}^{\gamma}}$. Then $\dot{H}^{0} = U$ and $\dot{H}^{\alpha} \subset \dot{H}^{\beta}$ if $\alpha \geq \beta$. Furthermore, we introduce the product space
\begin{align}\label{eq:product.space}
H^{\gamma} := \dot{H}^{\gamma} \times \dot{H}^{\gamma-1}, \quad \gamma \in \R,
\end{align}
endowed with the inner product
\begin{equation}
\langle X_1, X_2 \rangle_{H^{\gamma}} = \langle u_1,   u_2 \rangle_{\dot{H}^{\gamma}} + \langle v_1, v_2 \rangle_{\dot{H}^{\gamma-1}}, \quad X_1 = (u_1, v_1)^T, \: X_2 = (u_2, v_2)^T
\end{equation}
and the corresponding norm
\begin{align}\label{eq:product.norm}
\| X \|_{H^{\gamma} }^2 = \langle X, X \rangle_{H^{\gamma}} =  \|u\|_{\gamma}^2 + \| v\|^2_{\gamma-1}, \quad \gamma \in \mathbb{R},\: X=(u,v)^T.
\end{align}
For the special case $\gamma =0$, we denote $H := H^0 = \dot{H}^0 \times \dot{H}^{-1}$ and $(H,\, \langle \cdot, \cdot\rangle_H, \|\cdot\|_H )$ is a separable Hilbert space. Throughout this work we regard $\Lambda$ as an operator from $\dot{H}^1$ to $\dot{H}^{-1}$, defined by $( \Lambda u)(\varphi ) = \langle \nabla u, \nabla \varphi \rangle$ for $u, \varphi \in \dot{H}^1$, and define
\begin{align}\label{eq:Domain.A}
  \mathcal{D} (A) = \bigg\{ X = (u, v)^T \in H: A X =  \Big[\! \begin{array}{c}
    v
    \\
    -\Lambda u
    \end{array}\!\Big] \in H = \dot{H}^0 \times \dot{H}^{-1} \bigg\} = H^1 = \dot{H}^1 \times \dot{H}^{0}.
\end{align}
In this setting one can rigorously check that $A$ is closed and densely defined in $H$ and the resolvent set of $A$ contains all non-zero real numbers and $\|(\lambda \mathcal{I} - A)^{-1}\|_{\mathcal{L}(H)} \leq \tfrac{1}{|\lambda|}$ for any $\lambda \in \R$ (see, e.g., \cite[Section 5.3]{Lindgren12} for more details). Here and below, by $\mathcal{I} $ we mean the identity operator in $H$.  As a consequence, the operator $A$ is an infinitesimal generator of a $C_0$-group $E(t) = e^{t A}, t \in \R$ on $H$. In order to see the exact form of $E(t)$, note that $X(t) = E(t)X_0$ is the solution of the deterministic linear equation
\begin{equation}
\dot{X} = A X; \quad X(0) = X_0 = (u_0, v_0)^T.
\end{equation}
We solve it by using an eigenfunction expansion:
\begin{equation}\nonumber
X(t) = e^{tA} X_0 = \sum_{i = 1}^{\infty} \exp\! \Big(
    t \Big[\! \begin{array}{cc}
    0 & I
    \\
    -\lambda_i & 0
    \end{array}\!\Big]
    \Big) \Big[ \!\begin{array}{c}
    \langle u_0, e_i \rangle_U e_i
    \\
    \langle v_0, e_i \rangle_U e_i
    \end{array}\!\Big]
    =
    \sum_{i = 1}^{\infty}
    \bigg[\! \begin{array}{cc}
    \cos( \sqrt{\lambda_i} t )
    &
    \tfrac{1}{ \sqrt{\lambda_i}}\sin( \sqrt{\lambda_i} t )
    \\
    - \sqrt{\lambda_i} \sin( \sqrt{\lambda_i} t )
    &
    \cos( \sqrt{\lambda_i} t )
    \end{array}\!\bigg]
    \left[ \!\begin{array}{c}
    \langle u_0, e_i \rangle_U e_i
    \\
    \langle v_0, e_i \rangle_U e_i
    \end{array}\!\right].
\end{equation}
Hence the first component of $X = (u, v)^T$ is given by
\begin{equation}
\begin{split}
u(t) = & \sum_{i = 1}^{\infty} \Big[ \cos( \sqrt{\lambda_i} t ) \langle u_0, e_i \rangle_U e_i + \tfrac{1}{\sqrt{\lambda_i}} \sin( \sqrt{\lambda_i} t ) \langle v_0, e_i \rangle_U e_i \Big]
\\ = &
\cos(t\Lambda^{\frac{1}{2}}) u_0 + \Lambda^{-\frac{1}{2}}\sin(t\Lambda^{\frac{1}{2}}) v_0,
\end{split}
\end{equation}
and the second component of $X = (u, v)^T$ by
\begin{equation}
v(t) = \dot{u}(t) = -\Lambda^{\frac{1}{2}}\sin( t \Lambda^{\frac12} )u_0 +  \cos(t\Lambda^{\frac{1}{2}}) v_0.
\end{equation}
Here we introduced the cosine and sine operators $\cos(t\Lambda^{\frac{1}{2}})$ and $\sin(t\Lambda^{\frac{1}{2}})$. Accordingly, the semigroup $E(t) = e^{t A}$ should explicitly take the form as
\begin{align}\label{eq:Et}
E(t) = e^{tA} = \bigg[\! \begin{array}{cc}
    C(t) & \Lambda^{-\frac{1}{2}} S(t)
    \\
    -\Lambda^{\frac{1}{2}}S(t) & C(t)
    \end{array} \!\bigg],
\end{align}
where we further write $C(t)\!=\! \cos(t\Lambda^{\frac{1}{2}})$ and
$S(t)\!=\! \sin(t\Lambda^{\frac{1}{2}})$ for brevity.
Obviously, the cosine and sine operators satisfy a trigonometric identity in the sense that $\|S(t)u\|_U^2 + \|C(t)u\|_U^2 = \|u\|_U^2$ for $u \in U$. Using the trigonometric identity gives
\begin{equation} \label{eq:Et.group}
\|E(t) \|_{\mathcal{L}(H)} \leq 1, \quad t \in \R.
\end{equation}
Also, we will frequently use the fact that  $S(t), C(t), t \in \R$ and $\Lambda^{\gamma}, \gamma \in \R$ commute in the following estimates.
In the above setting, the scheme \eqref{eq:exp.Euler} can be rewritten as a recurrence equation in $H$:
\begin{equation}\label{eq:Abstr.Exp.Euler}
  X_{m+1} = E(\tau) \Big( X_m + \tau \BF (X_m) +  \BG (X_m) \Delta W_m \Big).
\end{equation}
For the convergence analysis, it is convenient to work with continuous processes.
Hence we define a continuous extension of \eqref{eq:Abstr.Exp.Euler}, $\tilde{X}(t) = (\tilde{u}(t), \tilde{v}(t))^T$, by
\begin{equation}\label{eq:X.tilde}
\tilde{X}(t) = E( t - t_m) \Big( X_m +  \BF (X_m) (t-t_m) + \BG (X_m) ( W(t) - W(t_m) ) \Big)
\end{equation}
for $ t \in [t_m, t_{m+1}]$.
It is obvious that $\tilde{X}(t)$ coincides with $X_m$ at the grid-points $t_m \in \mathcal{T}_M$. Note that \eqref{eq:Abstr.Exp.Euler} implies
\begin{equation}\label{eq:Xm}
\begin{split}
X_{m} = E(t_m) X_0 + \tau \sum_{k=0}^{m-1} E(t_m-t_k) \BF(X_k) + \sum_{k=0}^{m-1} E(t_m-t_k) \BG(X_k) \Delta W_k.
\end{split}
\end{equation}
Similarly, for arbitrary $t \in [0, T]$, $\tilde{X}(t) = (\tilde{u}(t), \tilde{v}(t))^T $ defined by \eqref{eq:X.tilde} can be expressed as
\begin{equation}\label{eq:X.tilde.global}
\tilde{X}(t) = E(t) X_0 + \int_0^t E(t- \lfloor s \rfloor_{\tau}) \BF ( X_{\lfloor s/\tau \rfloor}) \, \mbox{d} s + \int_{0}^t  E(t- \lfloor s \rfloor_{\tau}) \BG ( X_{\lfloor s/\tau \rfloor}) \, \mbox{d} W(s),
\end{equation}
where we define $\lfloor s/\tau \rfloor$ as an integer number not bigger than  $s/\tau$ and $\lfloor s \rfloor_{\tau} = \lfloor s/\tau \rfloor \cdot \tau$.


Next, we show the result on the existence and uniqueness of the mild solution of \eqref{eq:SEE}.

\begin{theorem}[Existence and uniqueness of mild solution]\label{thm:unique.mild}
Suppose all conditions in Assumption \ref{assump:FG} are fulfilled, let $W(t), t\in [0, T]$ be a cylindrical $Q$-Wiener process on the stochastic basis $\left(\Omega,\mathcal {F},\mathbb{P}, \{\mathcal {F}_t\}_{0\leq t\leq T}\right)$, given by \eqref{W.representation}, and let $X_0 = (u_0, v_0)^T$ be an $\mathcal{F}_0$-measurable $H$-valued random variable such that $\|X_0\|_{L^p(\Omega, H)} < \infty$ for some $ p \in [2, \infty)$. Then \eqref{eq:SEE} has a unique mild solution
\begin{align}\label{eq:mild.solution}
  X(t) = E(t)X_0 +  \int_0^t E(t-s) \BF(X(s)) \, \mbox{d}s + \int_0^t E(t-s)\BG(X(s)) \, \mbox{d} W(s), \: \: a.s.
\end{align}
for $t \in [0,T]$.
Moreover there exists a constant $C_{p,T} \in [0, \infty)$ depending on $p, T$ such that
\begin{align}\label{eq:X.MB}
  \sup_{t \in [0,T]} \big\|X(t)\big\|_{L^p(\Omega, H)} \leq C_{p,T} \big( \|X_0\|_{L^p(\Omega, H)} + 1 \big).
\end{align}
\end{theorem}
{\it Proof of Theorem \ref{thm:unique.mild}.}
Owing to \eqref{eq:Nf.condition1}, \eqref{eq:Nf.condition2} and
the definition of $\|\cdot \|_H$,
we infer that
\begin{align}
\| \BF(X_1)-\BF(X_2) \|_H =& \|\Lambda^{-\frac{1}{2}} \big(F(u_1)- F(u_2)\big)\|_U
\nonumber \\ \leq&
\|\Lambda^{-\frac{1}{2}}  \|_{\mathcal{L}(U)} \cdot L \|u_1-u_2\|_U \leq \tilde{L}\|X_1-X_2\|_H,
\label{eq:BF.LL} \\
\| \BF(X) \|_H =& \|\Lambda^{-\frac{1}{2}} F(u)\|_U \leq \tilde{L} (\|u\|_U + 1) \leq \tilde{L} (\|X\|_H+1)
\label{eq:BF.LG}
\end{align}
for arbitrary $X=(u, v)^T, X_1=(u_1, v_1)^T,
X_2=(u_2, v_2)^T$. Here $\tilde{L} = \|\Lambda^{-\frac{1}{2}}  \|_{\mathcal{L}(U)} \cdot L = L/\sqrt{\lambda_1}$.
Similarly, one can obtain that
\begin{align}
\| \BG(X_1)-\BG(X_2) \|_{\mathcal{L}_2(U_0, H)} = & \|\Lambda^{-\frac{1}{2}} \big(G(u_1)- G(u_2)\big)\|_{\mathcal{L}_2(U_0,U)}
\nonumber \\
\leq &
\| \Lambda^{\frac{-\delta}{2}} \|_{\mathcal{L}(U)} \cdot \|\Lambda^{\frac{\delta-1}{2}} \big(G(u_1)- G(u_2)\big)\|_{\mathcal{L}_2(U_0, U)}
\nonumber \\ \leq &
\hat{L} \|u_1-u_2\|_U \leq \hat{L} \|X_1-X_2\|_H,
\label{eq:BG.LL} \\
\| \BG(X) \|_{\mathcal{L}_2(U_0, H)} = & \|\Lambda^{-\frac{1}{2}} G(u)\|_{\mathcal{L}_2(U_0, U)} \leq \hat{L} (\|u\|_U + 1) \leq \hat{L} (\|X\|_H+1)
\label{eq:BG.LG}
\end{align}
with $\hat{L} =  L /\lambda_1^{\delta/2}$.
In view of Proposition 7.1 and Theorem 7.4 in \cite{DZ92}, the existence and uniqueness of the mild solution \eqref{eq:mild.solution} follow straightforwardly and \eqref{eq:X.MB} holds.
$\square$

In component-wise manner, the mild solution \eqref{eq:mild.solution} takes the form
\begin{align}\label{eq:mild.concrete}
\left\{\!
    \begin{array}{l}
    u(t)= C(t) u_0 + \Lambda^{-\frac{1}{2}} S(t)v_0 + \int_0^t \Lambda^{-\frac{1}{2}} S(t-s)F(u(s)) \,\mbox{d}s
    \\ \qquad \quad + \int_0^t \Lambda^{-\frac{1}{2}} S(t-s)G(u(s)) \,\mbox{d}W(s),
    \\
    v(t)= -\Lambda^{\frac{1}{2}} S(t) u_0 + C(t)v_0 + \int_0^t C(t-s)F(u(s)) \mbox{d}s
    \\ \qquad \quad
    + \int_0^t C(t-s)G(u(s)) \,\mbox{d}W(s)
    \end{array}\right.
\end{align}
$a.s.$ for $t \in [0, T]$. Therefore, $u(t)$ in \eqref{eq:mild.concrete} serves as the unique mild solution of \eqref{eq:Abstr.SWE}.


\section{Strong convergence}
\label{sect:stongconvergence}

This section focuses on the strong convergence of the scheme \eqref{eq:exp.Euler}. We begin by presenting a slightly modified version of the Burkholder-Davis-Gundy type inequality (\cite[Lemma 7.2]{DZ92}).

\begin{lemma} \label{lemma:BDG.inequality}
Let $(V,\, \langle \cdot, \cdot \rangle_V, \, \|\cdot \|_V)$ be a separable Hilbert space and let $ \Psi \colon [0, T] \times \Omega \rightarrow \mathcal{L}_2(U_0, V)$ be a predictable stochastic process.  Then for $t \in [0, T]$ and $p \in [2, \infty)$ there exists a constant $C_p$ such that
\begin{equation} \label{eq:BDG.inequality}
\bigg\| \sup_{s \in [0, t]} \Big\|\smallint_0^s \Psi(r) \, \mbox{d} W(r)\Big\|_{ V } \bigg\|_{ L^p(\Omega, \R) } \leq C_p \bigg( \smallint_0^t \Big\| \|\Psi(r)\|_{\mathcal{L}_2(U_0, V)} \Big\|^2_{L^p(\Omega, \R)} \mbox{d} r \bigg)^{\frac{1}{2}}.
\end{equation}
\end{lemma}
The following result is on further spatial regularity of the mild solution $X(t)$ and its numerical approximations.
\begin{prop}\label{prop:MB.X.Xm}
Assume all conditions in Theorem \ref{thm:unique.mild} are fulfilled and $X_0 \in L^p( \Omega, H^{\alpha})$ for some $\alpha \in [0, \infty)$ and $p \in [2, \infty)$. Then it holds for $\gamma \in [0,  \min(\alpha, \delta,  1)]$ that
\begin{equation}\label{eq:MB.Spatial.X}
\Big\| \sup_{ t \in[0, T]} \| X(t) \|_{ H^{\gamma}} \Big\|_{L^p(\Omega, \R)} \leq C (\|X_0\|_{L^p(\Omega, H^{\gamma})} + 1),
\end{equation}
and
\begin{equation}\label{eq:MB.Spatial.X.numer}
\Big\| \sup_{t \in [0,T]} \| \tilde{X}(t) \|_{ H^{\gamma}} \Big\|_{L^p(\Omega, \R)} \leq C (\|X_0\|_{L^p(\Omega, H^{\gamma} )} + 1).
\end{equation}
\end{prop}
{\it Proof of Proposition \ref{prop:MB.X.Xm}.} For the first step, we prove  \eqref{eq:MB.Spatial.X}. From \eqref{eq:mild.solution} we deduce  that
\begin{equation*}
\begin{split}
\sup_{t \in [0, T]} \|X(t)\|_{H^{\gamma} }
\leq &
\sup_{t \in [0, T]} \|E(t)X_0\|_{H^{\gamma} }
\\ &+
\sup_{t \in [0, T]} \|E(t- T)\|_{\mathcal{L}(H)} \cdot \sup_{t \in [0, T]} \Big\|\smallint_0^t E(T-s) \BF( X(s) ) \mbox{d}s \Big\|_{H^{\gamma}}
\\
& +
\sup_{t \in [0, T]} \|E(t- T)\|_{\mathcal{L}(H)} \cdot \sup_{t \in [0, T]} \Big\|\smallint_0^t E(T-s)\BG( X(s) ) \, \mbox{d} W(s)\Big\|_{ H^{\gamma} }.
\end{split}
\end{equation*}
Using \eqref{eq:Et.group} and elementary inequalities therefore results in
\begin{equation}
\begin{split}
\label{eq:X.MB.proof}
\Big\| \sup_{t \in [0, T]} \|X(t)\|_{H^{\gamma} } \Big\|_{L^p(\Omega, \R)}
\leq &
\|X_0\|_{L^p(\Omega, H^{\gamma}) } +  \smallint_0^T \left\| E(T-s) \BF( X(s) )\right\|_{ L^p(\Omega, H^{\gamma}) } \mbox{d}s
\\ & +
\Big\| \sup_{t \in [0, T]} \Big\|\smallint_0^t E(T-s)\BG( X(s) ) \, \mbox{d} W(s)\Big\|_{ H^{\gamma} } \Big\|_{ L^p(\Omega, \R) }
 \\ = &
\|X_0\|_{L^p(\Omega, H^{\gamma}) } + I_1 + I_2.
\end{split}
\end{equation}
Subsequently we will estimate $I_1$ and $I_2$ separately.
First, \eqref{eq:Nf.condition1}, \eqref{eq:X.MB}, the fact that $\gamma \leq \min(\alpha, \delta, 1)$, the trigonometric identity and elementary inequalities enable us to get
\begin{equation}
\begin{split}
I_1 = &
\smallint_0^T \big\| \big( \big\| \Lambda^{-\frac{1}{2}} S(T-s)F(u(s))\big\|^2_{\gamma} + \| C(T-s)F(u(s))\|^2_{\gamma-1}\big)^{1/2}  \big\|_{L^p(\Omega, \mathbb{R})} \mbox{d}s
\\ = &
\smallint_0^T \Big\|\big\| \Lambda^{\frac{\gamma - 1} {2}} F(u(s))\big\|_U \Big\|_{L^p(\Omega, \mathbb{R})} \mbox{d}s
\\ \leq &
\lambda_1^{\frac{\gamma-1}{2}} L \smallint_0^T \left( \|u(s) \|_{L^p(\Omega, U)} + 1 \right) \mbox{d}s
\\ \leq &
C (\|X_0\|_{L^p(\Omega, H)} + 1).
\end{split}
\end{equation}
Next, we note that
\begin{equation*}
\begin{split}
& \sup_{t \in [0, T]} \Big\|\smallint_0^t E(T-s)\BG( X(s) ) \, \mbox{d} W(s)\Big\|_{ H^{\gamma} }
\\ = &
\sup_{t \in [0, T]} \! \Big( \Big\| \smallint_0^t \Lambda^{-\frac{1}{2}} S(T-s) G(u(s)) \, \mbox{d}W(s) \Big\|^2_{\gamma} + \Big\| \smallint_0^t  C(T-s) G(u(s)) \, \mbox{d}W(s) \Big\|^2_{\gamma-1} \Big)^{\frac{1}{2}}
\\ \leq &
\sup_{t \in [0, T]} \! \Big\| \smallint_0^t \Lambda^{\frac{\gamma -1}{2}} S(T-s) G(u(s)) \, \mbox{d}W(s) \Big\|
+
\sup_{t \in [0, T]} \! \Big\| \smallint_0^t  \Lambda^{\frac{\gamma -1}{2}}  C(T-s) G(u(s)) \, \mbox{d}W(s) \Big\|.
\end{split}
\end{equation*}
Finally, we use the Burkholder-Davis-Gundy type inequality
\eqref{eq:BDG.inequality} and \eqref{eq:LHS} to arrive at
\begin{equation}
\begin{split}
|I_2|^2 \leq &
 2 C_p^2 \smallint_0^T \left\| \big\| \Lambda^{\frac{\gamma-1}{2}} S(T-s) G(u(s)) \big\|_{ \mathcal{L}_2^0} \right\|_{L^p(\Omega, \R)} ^2 \mbox{d}s
\\ & +
2 C_p^2 \smallint_0^T \left\| \big\| \Lambda^{\frac{\gamma-1}{2}} C(T-s) G(u(s)) \big\|_{ \mathcal{L}_2^0} \right\|_{L^p(\Omega, \R)} ^2 \mbox{d}s
\\  \leq &
4 C_p^2 \smallint_0^T \Big \| \big\| \Lambda^{\frac{\gamma-1}{2}} G(u(s)) \big\|_{ \mathcal{L}_2^0} \Big\|_{L^p(\Omega, \R)} ^2 \mbox{d}s
\\ \leq &
4 C_p^2 \smallint_0^T \big \| \Lambda^{\frac{\gamma - \delta}{2}} \big\|^2_{\mathcal{L}(U)} \cdot \Big \| \big\| \Lambda^{\frac{\delta-1}{2}} G(u(s)) \big\|_{ \mathcal{L}_2^0} \Big\|_{L^p(\Omega, \R)} ^2 \mbox{d}s
\\ \leq &
4 \lambda_1^{\gamma-\delta} C_p^2  \smallint_0^T \big \| L( \|u(s)\|_U + 1 ) \big\|_{L^p(\Omega, \R)} ^2 \mbox{d}s
\\
\leq &
8 L^2 \lambda_1^{\gamma-\delta} C_p^2 \smallint_0^T
\big( \| u(s) \|_{L^p(\Omega, U)} ^2 + 1 \big) \,
\mbox{d}s
\\
\leq &
C (\|X_0\|^2_{L^p(\Omega, H)} + 1).
\end{split}
\end{equation}
Here the stability properties of $S(t), C(t), \Lambda^{\frac{\gamma - \delta}{2}}$ for $t \in [0, T] $, $\gamma \leq \min( \alpha, \delta, 1)$ and some elementary inequalities were also used. Inserting the above estimates of $I_1$ and $I_2$ into \eqref{eq:X.MB.proof} and taking the assumption $X_0 \in L^p(\Omega, H^{\alpha})$ into account yield \eqref{eq:MB.Spatial.X}. To obtain \eqref{eq:MB.Spatial.X.numer}, we first derive from \eqref{eq:Xm} that
\begin{align}
\|X_m\|^2_{L^p(\Omega, H)} \leq & 3 \|E(t_m)X_0\|^2_{L^p(\Omega, H)} + 3m \tau^2 \sum_{k=0}^{m-1} \|E\big( t_m-t_k \big) \BF(X_k)\|^2_{L^p(\Omega, H)}
\nonumber \\ &
+ 3 C_p^2 \sum_{k=0}^{m-1} \smallint_{t_k}^{t_{k+1}} \big\| \| E\big( t_m-t_k  \big) \BG(X_k) \|_{L_2^0} \big\|^2_{L^p(\Omega, \R)} \mbox{d} s
\nonumber \\ \leq &
3 \|X_0\|^2_{L^p(\Omega, H)} + 3 T \tau \sum_{k=0}^{m-1}  \| \BF(X_k)\|^2_{L^p(\Omega, H)}
+ 3 C_p^2 \tau \sum_{k=0}^{m-1} \big\| \| \BG(X_k) \|_{L_2^0} \big\|^2_{L^p(\Omega, \R)}
\nonumber \\ \leq &
\hat{C}_{p,T} \big(\|X_0\|^2_{L^p(\Omega, H)} + 1\big) + \bar{C}_{p, T} \tau \sum_{k=0}^{m-1} \|X_k\|^2_{L^p(\Omega, H)},
\label{eq:Xm.MB}
\end{align}
where \eqref{eq:Et.group}, \eqref{eq:BF.LG} and \eqref{eq:BG.LG} were employed.
The discrete version of Gronwall's inequality applied to \eqref{eq:Xm.MB} and taking square roots show for all $m \in \{0, 1, \cdots, M\}$ that
\begin{equation}
\|X_m\|_{L^p(\Omega, H)} \leq C_{p,T} \big( \|X_0\|_{L^p(\Omega, H)} + 1 \big).
\end{equation}
With this and \eqref{eq:X.tilde.global}, the proof of  \eqref{eq:MB.Spatial.X.numer} goes along the same way as before.
$\square$

Assumption \ref{assump:FG} and Proposition \ref{prop:MB.X.Xm} together immediately imply two corollaries as follows.
\begin{cor} \label{cor:MB.FXGX}
Assume all conditions in Theorem \ref{thm:unique.mild} are fulfilled and $X_0 \in L^p(\Omega, H)$ for $p \in [2, \infty)$. Let $u(t)$ be given by \eqref{eq:mild.concrete}. Then there exists a constant $C$, depending on $p, L, T, \|X_0\|_{L^p(\Omega, H)}$, such that
\begin{align}\label{eq:MB.FXFXN}
\sup_{0\leq t \leq T} \big\| F (u(t)) \big\|_{L^p(\Omega, U)} \leq C, \quad \mbox{and} \quad  \sup_{0\leq t \leq T} \big\| \Lambda^{\frac{\delta-1}{2}} G (u(t)) \big\|_{L^p(\Omega, \mathcal{L}_2^0 )} \leq C.
\end{align}
\end{cor}

\begin{cor} \label{cor:MB.FXmGXm}
Assume all conditions in Theorem \ref{thm:unique.mild} are fulfilled and $X_0 \in L^p(\Omega, H)$ for $p \in [2, \infty)$. Let $u_m$ be given by \eqref{eq:exp.Euler}. Then there exists a constant $C$,  depending on $p, L, T, \|X_0\|_{L^p(\Omega, H)}$, such that
\begin{align}\label{eq:MB.FXmFXm}
\sup_{ 0 \leq m \leq M} \big \| F (u_m) \big\|_{L^p(\Omega, U)} \leq C, \quad \mbox{and} \quad  \sup_{ 0 \leq m \leq M} \big\| \Lambda^{\frac{\delta-1}{2}} G (u_m) \big\|_{L^p(\Omega, \mathcal{L}_2^0 )} \leq C.
\end{align}
\end{cor}
In what follows we present a useful lemma, which can be found in \cite{CLS13,KLL12}. The following H\"{o}lder continuity of $E(t)$ will put the ultimate limit 1 on the convergence order one can achieve for the time-stepping scheme \eqref{eq:exp.Euler}.
\begin{lemma} \label{lemma:SCR}
Assume that $S(t)$ and $C(t)$ are the sine and cosine operators and $E(t)$ is the group as defined above. Then for all $\gamma \in [0,1]$ there exists some constant $c_{\gamma}$ such that
\begin{align}\label{eq:SCNts}
\big\| \big( S(t)- S(s) \big) \Lambda^{-\frac{\gamma}{2}} \big\|_{\mathcal{L}(U)} \leq c_{\gamma}(t-s)^{\gamma},
\;
\big\| \big( C(t)- C(s) \big) \Lambda^{-\frac{\gamma}{2}} \big\|_{\mathcal{L}(U)} \leq c_{\gamma}(t-s)^{\gamma}
\end{align}
and
\begin{equation}\label{eq:Et.diff}
\| \big( E(t) - E(s)\big) X \|_H \leq c_{\gamma} (t-s)^{\gamma} \|X\|_{H^{\gamma}}
\end{equation}
for all $t \geq s \geq 0$.
\end{lemma}
Equipped with this lemma, one can investigate the H\"{o}lder regularity in time of the mild solution \eqref{eq:mild.concrete}, which is crucial in analyzing the approximation error of the time-discretization.

\begin{lemma}\label{lemma:ms.continuous}
Assume conditions in Theorem \ref{thm:unique.mild} are all fulfilled and $X_0 \in L^p( \Omega, H^{\alpha})$ for some $\alpha \in [0, \infty)$ and $p \in [2, \infty)$. Then it holds for $0 \leq s \leq t \leq T$ that
\begin{equation}\label{eq:u.holder.contin.}
\begin{split}
\left\|u(t)- u(s)\right\|_{L^p(\Omega, U)} \leq &
C(t-s)^{\min( \alpha, \delta, 1 )}.
\end{split}
\end{equation}
\end{lemma}
{\it Proof of Lemma \ref{lemma:ms.continuous}.} To get \eqref{eq:u.holder.contin.}, we first write
\begin{equation*}
X(t) - X(s) = \big( E(t-s) - \mathcal{I} \big) X(s) + \smallint_s^t E(t -r) \BF(X(r)) \, \mbox{d}r + \smallint_s^t E(t -r) \BG(X(r)) \, \mbox{d} W(r),
\end{equation*}
which admits
\begin{equation}
\begin{split}
 u(t)- u(s)
=&
\big(C(t-s)-I\big) u(s) + \Lambda^{-\frac{1}{2}} S(t-s)v(s)
\\ & +
\smallint_s^t \Lambda^{-\frac{1}{2}} S(t-r)F(u(r)) \, \mbox{d}r
+
 \smallint_s^t \Lambda^{-\frac{1}{2}} S(t-r)G(u(r)) \, \mbox{d}W(r).
\end{split}
\end{equation}
In the course of the proof, we assign $\rho = \min(\alpha, \delta, 1) $ for simplicity of presentation. The Burkholder-Davis-Gundy type inequality applied to the previous equality shows that
\begin{equation}\label{eq:ut-us2}
\begin{split}
\left\|u(t)- u(s)\right\|_{L^p(\Omega, U)} \leq &
\left\| \big(C(t-s)-I\big) u(s) \right\|_{L^p(\Omega, U)}
+
\big\|\Lambda^{-\frac{1}{2}} S(t-s) v(s) \big\|_{L^p(\Omega, U)}
\\ & +
\smallint_s^t \big\| \Lambda^{-\frac{1}{2}} S(t-r) F(u(r))\big\|_{L^p(\Omega, U)} \mbox{d}r
\\ & +
C_p \Big( \smallint_s^t \left\| \big\| \Lambda^{-\frac{1}{2}} S(t-r)G(u(r)) \big\|_{\mathcal{L}_2^0}\right\|^2_{L^p(\Omega, \R)} \mbox{d}r \Big)^{\frac{1}{2}}
\\ \leq &
\big\| \big(C(t-s)-I\big) \Lambda^{-\frac{\rho}{2}} \big\|_{\mathcal{L}(U)} \cdot \| u(s) \|_{L^p(\Omega, \dot{H}^{\rho})}
\\ & +
 \big\|\Lambda^{-\frac{\rho}{2}}S(t-s) \big\|_{\mathcal{L}(U)} \cdot \big\| v(s) \big\|_{L^p(\Omega, \dot{H}^{\rho-1})}
\\ & +
\smallint_s^t \! \big\| \Lambda^{-\frac{1}{2}} S(t-r) \big\|_{\mathcal{L}(U)} \cdot \big\| F(u(r))\big\|_{L^p(\Omega, U)} \mbox{d}r
\\ & +
C_p \Big( \smallint_s^t  \big\| \Lambda^{-\frac{\delta}{2}} S(t-r) \big\|^2_{\mathcal{L}(U)} \cdot \big\| \Lambda^{\frac{\delta -1}{2}} G(u(r)) \big\|^2_{ L^p( \Omega, \mathcal{L}_2^0) }  \mbox{d}r \Big)^{\frac{1}{2}}
\end{split}
\end{equation}
for $0 \leq s \leq t \leq T$ and $p \in [2, \infty)$. Note that Proposition \ref{prop:MB.X.Xm} ensures $\|u(s)\|_{L^p(\Omega, \dot{H}^{\rho})} \leq \|X(s)\|_{L^p(\Omega, H^{\rho})} < \infty$ and $\|v(s)\|_{L^p(\Omega, \dot{H}^{\rho-1})} \leq \|X(s)\|_{L^p(\Omega, H^{\rho})} < \infty$ for any $s \in [0, T]$.
Combining these bounds with Corollary \ref{cor:MB.FXGX} and Lemma \ref{lemma:SCR}, one can easily deduce from \eqref{eq:ut-us2} that
\begin{equation}
\begin{split}
\left\|u(t)- u(s)\right\|_{L^p(\Omega, U)} \leq &
C |t -s |^{\rho}
+
C \smallint_s^t  (t-r)\, \| F(u(r))\|_{L^p(\Omega, U)} \mbox{d}r
\\& +
C \Big( \smallint_s^t  (t -r)^{2\min( \delta, 1) } \big\| \Lambda^{\frac{\delta -1}{2}} G(u(r)) \big\|^2_{L^p( \Omega, \mathcal{L}_2^0) }  \mbox{d}r \Big)^{\frac{1}{2}}
\\ \leq &
C |t -s |^{\rho}
\end{split}
\end{equation}
for $0 \leq s \leq t \leq T$ and $p \in [2, \infty)$.
This finishes the proof of Lemma \ref{lemma:ms.continuous}. $\square$

It is worthwhile to remark that one can similarly work with $X(t)$ and examine its H\"{o}lder regularity in the product space $H$. But this leads to reduced H\"{o}lder regularity in time and one can only obtain $\|X(t) - X(s)\|_{L^p(\Omega, H)} \leq C (t-s)^{\min( \alpha, \delta, \frac{1}{2})}$, which puts an ultimate limit $\tfrac{1}{2}$ on the strong convergence order of the scheme \eqref{eq:exp.Euler}.
Now we formulate the main result in this section as follows.
\begin{theorem}\label{thm:strong.conv.result}
Suppose that all conditions in Theorem \ref{thm:unique.mild} are fulfilled and assume $X_0 \in L^p(\Omega, H^1)$ for some $p \in [2, \infty)$.
Then it holds that
\begin{equation} \label{eq:strong.conv.result}
\Big\| \sup_{t \in [0, T] } \big\| X(t) - \tilde{X}(t) \big\|_H \Big\|_{L^p(\Omega, \R)} \leq C \tau^{\min( \delta, 1)},
\end{equation}
where $X(t)$ is the mild solution of \eqref{eq:SEE} and $\tilde{X}(t)$ is a continuous-time extension of $X_m$, given by \eqref{eq:X.tilde}.
\end{theorem}

As an immediate consequence, we have the following strong convergence result.

\begin{cor}\label{cor:strong.conv.result}
Suppose that all conditions in Theorem \ref{thm:strong.conv.result} are fulfilled.
Then
\begin{equation} \label{eq:cor.strong.conv.result}
\Big\| \sup_{t \in [0, T] } \big\| u(t) - \tilde{u}(t) \big\|_U \Big\|_{L^p(\Omega, \R)} \leq C \tau^{\min( \delta, 1)},
\end{equation}
where $u(t)$ and $\tilde{u}(t)$ are the first components of $X(t)$ and $\tilde{X}(t)$, respectively.
\end{cor}
The above results reveal that the order of strong convergence is essentially governed by the spatial regularity of the noise term and is in accordance with the exponents of the H\"{o}lder regularity in time as stated in Lemma \ref{lemma:ms.continuous}. In addition, one can easily observe that the upper limit on strong order can be achieved only when $\delta \in [1, \infty)$. Furthermore, as indicated in \cite{CLS13} for the linear SWE, the exponential scheme \eqref{eq:exp.Euler} for nonlinear SWE also allows for higher strong order than the implicit Euler and Crank-Nicolson time discretizations do (compare Corollary \ref{cor:strong.conv.result} with \cite[Theorem 4.6]{KLL12}).

{\it Proof of Theorem \ref{thm:strong.conv.result}. }
Combining \eqref{eq:X.tilde.global} and \eqref{eq:mild.solution} gives
\begin{equation}
\begin{split}
\label{eq:X.time.diff2}
 X(s) - \tilde{X}(s)  = &
\smallint_{0}^{s } \Big( E(s -r) \BF(X(r)) - E(s - \lfloor r \rfloor_{\tau}) \BF(X_{ \lfloor r/\tau \rfloor }) \Big) \,\mbox{d} r
 \\  & +
 \smallint_{0}^{s } \Big( E(s -r) \BG(X(r)) - E(s - \lfloor r \rfloor_{\tau}) \BG(X_{ \lfloor r/\tau \rfloor }) \Big) \,\mbox{d} W(r).
\end{split}
\end{equation}
Thanks to the same arguments as used in  \eqref{eq:X.MB.proof}, one can similarly get
\begin{align}
\label{eq:varepsilon.m}
& \Big\| \sup_{s \in [0, t] } \big\| X(s) - \tilde{X}(s) \big\|_H \Big\|_{L^p(\Omega, \R)}
\nonumber \\ \leq &
\smallint_{0}^{t}  \Big \|  E(t-r) \BF(X(r)) - E(t- \lfloor r \rfloor_{\tau} ) \BF(X_{\lfloor r/\tau \rfloor} ) \Big\|_{L^p(\Omega, H)}\, \mbox{d} r
\nonumber \\ & +
\bigg \| \sup_{s \in [0, t]} \Big\| \smallint_{0}^{s} \Big( E(t-r) \BG(X(r)) - E(t- \lfloor r \rfloor_{\tau} ) \BG(X_{\lfloor r/\tau \rfloor}) \Big) \, \mbox{d} W(r) \Big\|_H \bigg\|_{L^p(\Omega, \R)}
\nonumber \\ = &
I_3
+  I_4.
\end{align}
Noting that $\tilde{X} (\lfloor r \rfloor_{\tau}) = X_{ \lfloor r/\tau \rfloor }$ and using \eqref{eq:Et.group}, \eqref{eq:BF.LL}, \eqref{eq:MB.FXFXN}, \eqref{eq:u.holder.contin.} and Lemma \ref{lemma:SCR}, we find that
\begin{equation}\label{eq:I_3}
\begin{split}
I_3 \leq &
\smallint_{0}^{t} \left\|  \Big( E(t-r)  - E(t-\lfloor r \rfloor_{\tau} ) \Big) \BF(X(r)) \right\|_{L^p(\Omega, H)}  \mbox{d} r
\\ & +
\smallint_{0}^{t} \left\| E(t- \lfloor r \rfloor_{\tau} ) \Big( \BF(X(r)) -  \BF(X( \lfloor r \rfloor_{\tau} )) \Big)\right\|_{L^p(\Omega, H)}  \mbox{d} r
\\ & +
\smallint_{0}^{t} \left\| E(t- \lfloor r \rfloor_{\tau} ) \Big( \BF(X( \lfloor r \rfloor_{\tau} )) - \BF( X_{\lfloor r/\tau \rfloor} ) \Big)\right\|_{L^p(\Omega, H)}  \mbox{d} r
\\ \leq &
c_1 \tau  \smallint_0^t \big\| \BF(X(r)) \big\|_{L^p(\Omega, H^1)} \mbox{d} r
+
\smallint_{0}^{t} \big\| \BF(X(r)) -  \BF(X( \lfloor r \rfloor_{\tau} )) \big\|_{L^p(\Omega, H)}  \mbox{d} r
\\ & +
\smallint_{0}^{t} \big\| \BF(X(\lfloor r \rfloor_{\tau})) -  \BF(\tilde{X}( \lfloor r \rfloor_{\tau} )) \big\|_{L^p(\Omega, H)}  \mbox{d} r
\\ \leq &
c_1 \tau \smallint_0^t \|F(u(r))\|_{L^p(\Omega, U)} \, \mbox{d} r
+
L/\sqrt{\lambda_1} \smallint_{0}^{t} \big\| u(r) -  u( \lfloor r \rfloor_{\tau} ) \big\|_{L^p(\Omega, U)}  \mbox{d} r
\\ & +
L/\sqrt{\lambda_1} \smallint_{0}^{t} \big\| X(\lfloor r \rfloor_{\tau}) -  \tilde{X}( \lfloor r \rfloor_{\tau} ) \big\|_{L^p(\Omega, H)}  \mbox{d} r
\\ \leq &
C \tau + C \tau^{\min(\delta, 1)} + C \smallint_{0}^{t} \Big\| \sup_{r \in [0, s] } \big\| X(r) - \tilde{X}(r) \big\|_H \Big\|_{L^p(\Omega, \R)}  \mbox{d} s
\\ \leq &
C \tau^{\min(\delta, 1)} + C \smallint_{0}^{t} \Big\| \sup_{r \in [0, s] } \big\| X(r) - \tilde{X}(r) \big\|_H \Big\|_{L^p(\Omega, \R)}  \mbox{d} s.
\end{split}
\end{equation}
With regard to $I_4$, we apply the Burkholder-Davis-Gundy type inequality \eqref{eq:BDG.inequality} to get
\begin{equation}
\begin{split}
|I_4|^2 \leq&  C_p^2 \smallint_{0}^{t} \Big \| E(t-r) \BG(X(r)) - E(t- \lfloor r \rfloor_{\tau} ) \BG(X_{\lfloor r/\tau \rfloor}) \Big\|^2_{L^p(\Omega, \mathcal{L}_2^0)} \, \mbox{d} r
\\ \leq &
3 C_p^2 \smallint_{0}^{t} \left\|  \Big( E(t-r)  - E(t-\lfloor r \rfloor_{\tau} ) \Big) \BG(X(r)) \right\|^2_{L^p(\Omega, \mathcal{L}_2^0)}  \mbox{d} r
\\ & +
3 C_p^2 \smallint_{0}^{t} \left\| E(t- \lfloor r \rfloor_{\tau} ) \Big( \BG(X(r)) -  \BG(X( \lfloor r \rfloor_{\tau} )) \Big)\right\|^2_{L^p(\Omega, \mathcal{L}_2^0)}  \mbox{d} r
\\ & +
3 C_p^2 \smallint_{0}^{t} \left\| E(t- \lfloor r \rfloor_{\tau} ) \Big( \BG(X( \lfloor r \rfloor_{\tau} )) - \BG( X_{\lfloor r/\tau \rfloor} ) \Big)\right\|^2_{L^p(\Omega, \mathcal{L}_2^0)}  \mbox{d} r.
\end{split}
\end{equation}
Further, employing \eqref{eq:Et.group}, \eqref{eq:BG.LL}, \eqref{eq:MB.FXFXN}, \eqref{eq:u.holder.contin.} and Lemma \ref{lemma:SCR} shows that
\begin{equation}\label{eq:I_4}
\begin{split}
|I_4|^2 \leq& c_{\delta} \tau^{2\min(\delta, 1)} \smallint_0^t \Big\|  \| \BG(X(r)) \|_{\mathcal{L}_2(U_0, H^{\delta})} \Big\|^2_{L^p(\Omega, \R)}  \mbox{d} r
\\ & +
3 C_p^2 \smallint_{0}^{t} \left\| \BG(X(r)) -  \BG(X( \lfloor r \rfloor_{\tau} )) \right\|^2_{L^p(\Omega, \mathcal{L}_2^0)}  \mbox{d} r
\\ & +
3 C_p^2 \smallint_{0}^{t} \left\| \BG(X( \lfloor r \rfloor_{\tau} )) - \BG( X_{\lfloor r/\tau \rfloor} )\right\|^2_{L^p(\Omega, \mathcal{L}_2^0)}  \mbox{d} r
\\ \leq&
c_{\delta} \tau^{2\min(\delta, 1)} \smallint_0^t \Big\|  \| \Lambda^{\frac{\delta-1}{2}} G(u(r)) \|_{\mathcal{L}_2(U_0, U)} \Big\|^2_{L^p(\Omega, \R)}  \mbox{d} r
\\ & +
3 C_p^2 L^2 /\lambda_1^{\delta} \smallint_{0}^{t} \left\| u(r) -  u( \lfloor r \rfloor_{\tau} ) \right\|^2_{L^p(\Omega, U)}  \mbox{d} r
\\ & +
3 C_p^2 L^2 /\lambda_1^{\delta} \smallint_{0}^{t} \big\| X( \lfloor r \rfloor_{\tau} ) - \tilde{X}( \lfloor r \rfloor_{\tau} )\big\|^2_{L^p(\Omega, H)}  \mbox{d} r
\\ \leq &
C \tau^{2\min(\delta, 1)} + C \smallint_{0}^{t} \Big\| \sup_{r \in [0, s] } \big\| X(r) - \tilde{X}(r) \big\|_H \Big\|^2_{L^p(\Omega, \R)}  \mbox{d} s.
\end{split}
\end{equation}
Defining a non-decreasing function $\vartheta \colon [0, T] \rightarrow \R $ by $\vartheta(t) = \big\| \sup_{s \in [0, t] } \| X(s) - \tilde{X}(s) \|_H \big\|^2_{L^p(\Omega, \R)}$ and taking the estimates \eqref{eq:I_3} and \eqref{eq:I_4} into account we derive from \eqref{eq:varepsilon.m} that
\begin{equation}
\begin{split}
\vartheta(t) = \Big\| \sup_{s \in [0, t] } \big\| X(s) - \tilde{X}(s) \big\|_H \Big\|^2_{L^p(\Omega, \R)}
\leq  C \tau^{2 \min(\delta, 1)} + C \int_0^t \vartheta(s)\, \mbox{d} s.
\end{split}
\end{equation}
Proposition \ref{prop:MB.X.Xm} guarantees the boundedness of $\vartheta(t) $ for $t \in [0, T]$, which enables us to apply Gronwall's lemma (see, e.g., \cite[Theorem 8.1]{Mao97}) to show for all $t \in [0, T]$ that
\begin{equation}\label{eq:mean.square.error}
\vartheta(t) = \Big\| \sup_{s \in [0, t] } \big\| X(s) - \tilde{X}(s) \big\|_H \Big\|^2_{L^p(\Omega, \R)} \leq C \tau^{2 \min( \delta, 1)}.
\end{equation}
Finally, taking square roots of \eqref{eq:mean.square.error} completes the proof of Theorem \ref{thm:strong.conv.result}. $\square$

\section{Applications to examples}
\label{sec:examples}

The aim of this section is to include several concrete examples which fit in the abstract setting formulated above. To this end, let $\mathcal{O} \subset \mathbb{R}^d$, $d = 1,2,3$, be a bounded open set with smooth boundary $\partial \mathcal{O}$. A nonlinear white noise driven stochastic wave equation (SWE) with Dirichlet boundary condition is usually described by
\begin{equation}\label{eq:SWE}
\left\{
    \begin{array}{lll}
    \frac{\partial^2 u}{\partial t^2} = \Delta u + f(\xi, u) + g(\xi, u)\dot{W}(t), \quad  t \in (0, T], \:\: \xi \in \mathcal{O},
    \\
     u(0, \xi) = u_0(\xi), \, \frac{\partial u}{\partial t}  (0,\xi) = v_0(\xi), \: \xi \in \mathcal{O},
     \\
     u|_{\partial \mathcal{O}} = 0,  \: t \in (0, T],
    \end{array}\right.
\end{equation}
where $T \in (0, \infty)$, $\Delta = \sum_{k=1}^d \frac{\partial^2}{\partial \xi_k^2}$ is the Laplace operator and $f,g \colon \mathcal{O} \times \mathbb{R} \rightarrow \mathbb{R}$ are deterministic functions. Moreover, the initial data $u_0, v_0 \colon \mathcal{O} \times \Omega \rightarrow \mathbb{R} $ are random variables and $W$ is a noise process, which will be specified later. Numerics of such equation has been considered in \cite{CL07,HE10,QS06,Walsh06,WGT13}.

Let $U := L^2\big(\mathcal{O}, \mathbb{R}\big)$ be the separable Hilbert space of real-valued square integrable functions from $\mathcal{O}$ to $\R$, with the scalar product and the norm
\begin{equation*}
\langle u_1, u_2 \rangle_U = \int_{\mathcal{O}} u_1(\xi ) u_2(\xi) \mbox{d} \xi, \quad \|u\|_U = \left( \int_{\mathcal{O}} |u(\xi )|^2 \mbox{d} \xi \right)^{1/2}
\end{equation*}
for all $u, u_1, u_2 \in U$. Moreover, let $\{W(t)\}_{t \in [0, T]}$ be a cylindrical $Q$-Wiener process on a stochastic basis $\left(\Omega,\mathcal {F},\mathbb{P}, \{\mathcal {F}_t\}_{t\in [0,T]}\right)$, given by \eqref{W.representation}. It is a classical result that the covariance operator $Q = Q^{\frac{1}{2}} \circ Q^{\frac{1}{2}} \in \mathcal{L}(U)$ is a nonnegative, symmetric operator so that
\begin{equation}
Q \phi_i = q_i \phi_i, \, q_i \geq 0, \, i \in \mathbb{N}.
\end{equation}
Then one can rewrite \eqref{eq:SWE} in an abstract It\^{o} form as \eqref{eq:Abstr.SWE}, namely,
\begin{align*}
\left\{\!
    \begin{array}{ll}
    \mbox{d} \dot{u} = -\Lambda u \, \mbox{d}t + F(u)\, \mbox{d}t + G(u)\,\mbox{d}W(t), \quad t \in (0, T],
    \\
     u(0) = u_0, \, \dot{u}(0) = v_0,
    \end{array}\right.
\end{align*}
where $-\Lambda : \mathcal{D}(\Lambda)\subset U \rightarrow U $ denotes the Laplacian with homogeneous Dirichlet boundary condition, and where $F\colon U \rightarrow U$, $G(u) \colon U \rightarrow U $ are the Nemytskij operators given by
\begin{equation}\label{eq:Nemytskij}
F (u)(\xi) = f(\xi, u(\xi)),  \quad \big(G (u)(\varphi) \big)(\xi) = g(\xi, u(\xi)) \cdot \varphi(\xi), \quad \xi \in \mathcal{O}.
\end{equation}

With the above setting, we take a close look at conditions in Assumption \ref{assump:FG}. To do this we consider separately several cases as follows.


\subsection{Additive noise}\label{subsec:additive.noise}

For the first case, we consider SWEs driven by additive noise with $g(\xi, u) \equiv 1 $ for all $\xi \in \mathcal{O}, u \in \R$. Furthermore, assume for some $\beta \in (0, \infty)$ that
\begin{equation}\label{eq:Lambda.betta}
\| \Lambda^{\frac{\beta -1}{2}} Q^{\frac{1}{2} } \|_{\mathcal{L}_2(U)} = \bigg( \sum_{i = 1}^{\infty} q_i \| \Lambda^{\frac{\beta -1}{2}} \phi_i \|_U^2 \bigg)^{\frac{1}{2}}  < \infty.
\end{equation}
Concerning $f$, we assume $f\colon \mathcal{O} \times \R \rightarrow \R $ in \eqref{eq:SWE} satisfies
\begin{equation}\label{eq:f.LL.GL}
\begin{split}
|f(\xi, u)| \leq  L(|u|+1), \quad
|f(\xi, u_1) - f(\xi, u_2)| \leq&  L |u_1 - u_2|
\end{split}
\end{equation}
for all $\xi \in \mathcal{O}$, $u, u_1, u_2 \in \R$.
Then the Nemytskij operator $F$ satisfies
\begin{equation}
\begin{split}
\| F(u) \|_U^2 =& \int_{\mathcal{O}} \big| f(\xi, u(\xi) ) \big|^2 \mbox{d} \xi
\leq
2L^2 \int_{\mathcal{O}} \big( |u(\xi)|^2 + 1 \big) \mbox{d} \xi = 2L^2 \big( \|u\|_U^2 + \varrho(\mathcal{O}) \big),
\end{split}
\end{equation}
where $\varrho(\mathcal{O}) $, the measure of the set $\mathcal{O}$, is bounded by assumption.
In a similar way,
\begin{equation}
\big\| F(u_1)-F(u_2) \big\|_U  \leq L \|u_1 - u_2\|_U.
\end{equation}
Moreover, $G(u) \equiv I$ in this setting and thus $G(u_1) - G(u_2) \equiv 0$ and
\begin{equation}
\big\| \Lambda^{\frac{\beta -1 }{2}} G(u)  \big\|_{\mathcal{L}_2^0} = \big\| \Lambda^{\frac{\beta -1 }{2}} \big\|_{\mathcal{L}_2^0} = \| \Lambda^{\frac{\beta -1}{2}} Q^{\frac{1}{2} } \|_{\mathcal{L}_2(U)} < \infty.
\end{equation}
Therefore, Assumption \ref{assump:FG} is fulfilled with $\delta = \beta$. An immediate consequence of Corollary \ref{cor:strong.conv.result} reads:
\begin{cor}\label{cor:str.conv.additive}
Suppose that $g(\xi, u) \equiv 1 $ for all $\xi \in \mathcal{O}$, $u \in \R$ and $f \colon \mathcal{O} \times \R \rightarrow \R$ satisfies \eqref{eq:f.LL.GL}. Let $W(t)$ be a cylindrical $Q$-Wiener process on the stochastic basis $\left(\Omega,\mathcal {F},\mathbb{P}, \{\mathcal {F}_t\}_{0\leq t\leq T}\right)$ with \eqref{eq:Lambda.betta} fulfilled. Additionally, we assume $u_0, v_0 \in \mathcal{F}_0$ and $ u_0 \in L^p(\Omega, \dot{H}^1), v_0 \in L^p(\Omega, U)$ for some $p \in [2, \infty)$. Then the problem \eqref{eq:SWE} has a unique mild solution. Moreover, there exists a generic constant $C \in [0, \infty)$ depending on $T, L, p, \beta$ and the initial data, such that
\begin{equation}\label{eq:strong.conv.additive}
\Big\| \sup_{t_m \in \mathcal{T}_M } \big\| u(t_m) - u_m \big\|_{U} \Big\|_{L^p(\Omega, \R)} \leq C \tau^{\min( \beta, 1) },
\end{equation}
where $u(t)$ is the mild solution of \eqref{eq:SWE} and $u_m$ is the numerical solution produced by \eqref{eq:exp.Euler}.
\end{cor}

\subsection{Multiplicative trace class noise }

Assume $f \colon \mathcal{O} \times \R \rightarrow \R $ in \eqref{eq:SWE} satisfies \eqref{eq:f.LL.GL} and  $g \colon \mathcal{O} \times \R \rightarrow \R $ in \eqref{eq:SWE} satisfies
\begin{equation}\label{eq:g.LL.GL}
\begin{split}
|g(\xi, u)| \leq  L(|u|+1), \quad
|g(\xi, u_1) - g(\xi, u_2)| \leq&  L |u_1 - u_2|
\end{split}
\end{equation}
for all $\xi \in \mathcal{O}$, $u, u_1, u_2 \in \R$.
Furthermore, we assume in this subsection
\begin{equation}\label{eq:Bounds.Qi}
\text{Tr} ( Q ) = \sum_{i \in \N} q_i < \infty, \quad \sup_{i \in \N}\sup_{ \xi \in \bar{\mathcal{O}} } |\phi_i( \xi )| < \infty.
\end{equation}
Then it is not difficult to see that
\begin{equation}
\begin{split}
\|G(u)\|^2_{\mathcal{L}_2^0}=& \big\| G(u)Q^{\frac{1}{2}} \big\|^2_{\mathcal{L}_2(U)} = \sum_{i=1}^{\infty} \big\| G(u) Q^{\frac{1}{2}} \phi_i\big\|^2_U
\\ =&
\sum_{i=1}^{\infty} q_i \big\| G(u) \phi_i\big\|^2_U
=
\sum_{i=1}^{\infty} q_i \int_{\mathcal{O}} \big| g ( \xi, u(\xi) ) \phi_i(\xi) \big|^2 \mbox{d} \xi
\\ \leq &
\sup_{i \in \N}\sup_{ \xi \in \bar{\mathcal{O}} } |\phi_i( \xi )|^2 \cdot \mbox{Tr} (Q) \cdot 2 L^2 \big( \|u\|_U^2 + \varrho(\mathcal{O}) \big).
\end{split}
\end{equation}
In the same way, one can obtain
\begin{equation}
\big\| \big( G(u_1)-G(u_2) \big) Q^{\frac{1}{2}} \big\|^2_{\mathcal{L}_2(U)}
\leq
\sup_{i \in \N}\sup_{ \xi \in \bar{\mathcal{O}} } |\phi_i( \xi )|^2 \cdot \mbox{Tr} (Q) \cdot L^2 \|u_1 - u_2\|_U^2.
\end{equation}
Consequently, Assumption \ref{assump:FG} is fulfilled with $\delta = 1$ and the following corollary follows.

\begin{cor}\label{cor:str.conv.trace}
Suppose that $f, g \colon \mathcal{O} \times \R \rightarrow \R$ satisfy \eqref{eq:f.LL.GL} and \eqref{eq:g.LL.GL}. Let $W(t)$ be a standard $U$-valued $Q$-Wiener process on a stochastic basis $\left(\Omega,\mathcal {F},\mathbb{P}, \{\mathcal {F}_t\}_{0\leq t\leq T}\right)$ with \eqref{eq:Bounds.Qi} fulfilled.  Additionally, we assume $u_0, v_0 \in \mathcal{F}_0$ and $ u_0 \in L^p(\Omega, \dot{H}^1), v_0 \in L^p(\Omega, U)$ for some $p \in [2, \infty)$. Then the problem \eqref{eq:SWE} has a unique mild solution and it holds that
\begin{equation}
\Big\| \sup_{t_m \in \mathcal{T}_M } \big\| u(t_m) - u_m \big\|_{U} \Big\|_{L^p(\Omega, \R)} \leq C \tau,
\end{equation}
where $u(t)$ is the mild solution of \eqref{eq:SWE} and $u_m$ is produced by \eqref{eq:exp.Euler}.
\end{cor}

\subsection{Multiplicative space-time white noise }

In this subsection, let $d = 1$, $\mathcal{O} = (0, 1)$ and let $Q = I$. Then $W(t)$ becomes a cylindrical $I$-Wiener process, which can be given by
\begin{equation}\label{eq:Wt.space.time.white}
W(t) = \sum_{i\in \N} \beta_i(t) e_i, \quad t \in [0, T],
\end{equation}
where $\{\beta_i(t) \}_{i\in \mathbb{N}}$ are a family of mutually independent real Brownian motions and $\{e_i=\sqrt{2}\sin(i \pi x), \,x \in(0,1)\}_{i \in \mathbb{N}}$ form an orthonormal basis of $U$ consisting of eigenfunctions of $\Lambda$ with $\Lambda e_i = \lambda_i e_i,\, \lambda_i=\pi^2 i^2, \, i \in \mathbb{N}$.
It is then obvious for any  $\epsilon >0$ that
\begin{align}\label{eq:Lambda.Q}
\big\|\Lambda^{-\frac{\epsilon+1}{4}} \big\|^2_{\mathcal{L}_2(U)} = \pi^{-(\epsilon+1)} \sum_{i=1}^{\infty} i^{-(\epsilon + 1)}  < \infty, \quad \mbox{and} \quad \sup_{i \in \N } \sup_{x \in [0, 1]} |e_i(x)| \leq \sqrt{2}.
\end{align}
Therefore, it follows that
\begin{equation}
\begin{split}
\big\| \Lambda^{-\frac{\epsilon+1}{4}} G(u) \big\|^2_{\mathcal{L}_2(U)} =&
\sum_{i=1}^{\infty} \big\| \Lambda^{-\frac{\epsilon+1}{4}} G(u) e_i\big\|^2_U
=
\sum_{i=1}^{\infty} \sum_{j = 1}^{\infty} \big| \langle \Lambda^{-\frac{\epsilon+1}{4}} G(u) e_i, e_j\rangle \big|^2
\\= &
\sum_{i=1}^{\infty} \sum_{j = 1}^{\infty} \big| \langle G(u) e_i, \lambda_j^{-\frac{\epsilon+1}{4}} e_j\rangle \big|^2
\\ = &
\sum_{i=1}^{\infty} \sum_{j = 1}^{\infty} \lambda_j^{-\frac{\epsilon+1}{2} } \Big| \int_{\mathcal{O}}  g ( \xi, u(\xi) ) e_i(\xi) e_j(\xi) \, \mbox{d} \xi \Big|^2
\\ \leq &
2 \sum_{i=1}^{\infty} \sum_{j = 1}^{\infty} \lambda_j^{-\frac{\epsilon+1}{2} } \Big| \int_{\mathcal{O}}  g ( \xi, u(\xi) ) e_i(\xi) \, \mbox{d} \xi \Big|^2
\\ = &
2 \pi^{-(\epsilon+1) } \sum_{j = 1}^{\infty} j^{-(\epsilon+1) } \cdot \|g(\cdot, u(\cdot))\|_U^2 \leq \hat{C}_{\epsilon} \big( \|u\|_U^2 + \varrho(\mathcal{O}) \big)
\end{split}
\end{equation}
for arbitrarily small $\epsilon >0$.
Similarly, one can show for arbitrarily small $\epsilon >0$ that
\begin{equation}
\big\| \Lambda^{-\frac{\epsilon+1}{4}} \big( G(u_1) - G(u_2) \big) \big\|_{\mathcal{L}_2(U)} \leq \check{C}_{\epsilon} \|u_1 - u_2 \|_{U}.
\end{equation}
Hence, Assumption \ref{assump:FG} is satisfied with $ \delta = \frac{1-\epsilon}{2}$ for arbitrarily small $\epsilon >0$.

\begin{cor}\label{cor:Strong.Converg.Spa.time.white}
Suppose that $f, g \colon (0, 1) \times \R \rightarrow \R$ satisfy \eqref{eq:f.LL.GL} and \eqref{eq:g.LL.GL}. Additionally, we assume $u_0, v_0 \in \mathcal{F}_0$ and $ u_0 \in L^p(\Omega, \dot{H}^1), v_0 \in L^p(\Omega, U)$ for some $p \in [2, \infty)$. Moreover, let $W(t)$ be a cylindrical $I$-Wiener process given by \eqref{eq:Wt.space.time.white}. Then the problem \eqref{eq:SWE} has a unique mild solution. Moreover, there exists a constant $C \in [0, \infty)$ depending on $T, L, p, \epsilon$ and the initial data, such that
\begin{equation}
\Big\| \sup_{t_m \in \mathcal{T}_M } \big\| u(t_m) - u_m \big\|_{U} \Big\|_{L^p(\Omega, \R)} \leq C_{\epsilon} \tau^{\frac{1- \epsilon}{2} },
\end{equation}
where $u(t)$ is the mild solution of \eqref{eq:SWE} and $u_m$ is the numerical solution produced by \eqref{eq:exp.Euler}.
\end{cor}

%
%

\section{Weak convergence}
\label{sec:weak.convergence}

In many applications, a key task is to approximate the quantity $\mathbb{E} \! \left[ \varphi (u(T))\right]$,
where $\varphi$ is a functional of the mild solution to \eqref{eq:Abstr.SWE}. This leads to another important notion of convergence for a numerical scheme,  the weak convergence, which is concerned with the approximation of law.
In this section, let us turn our attention to this topic. As is well known that the rate of the weak error can, in some situations, be improved compared to that of the strong error. Below, we shall show this fact for nonlinear SWE driven by additive noise. Since the treatment of multiplicative noise case causes additional technical difficulties, it will be addressed elsewhere in our future work.

For the additive case when $G(u) \equiv I$ for $u \in U$, the abstract equation \eqref{eq:SEE} reduces to
\begin{equation}\label{eq:SEE.additive}
\left\{
    \begin{array}{ll}
    \mbox{d} X(t) = A X(t)\, \mbox{d}t + \mathbf{F}(X)\, \mbox{d}t + B\, \mbox{d}W(t), \quad  t \in (0, T],
    \\
     X(0)= X_0,
    \end{array}\right.
\end{equation}
where we adopt the same notations as in \eqref{eq:SEE.notations} and additionally define here $B \colon U \rightarrow H$
by $B = (0, I)^T$.

To begin with, we need the following assumption.
\begin{assumption}\label{ass:weak.error.assump}
Assume $G(u) \equiv I$ for all $u \in U$, and let $F \colon U \rightarrow U $ be a twice differentiable mapping satisfying
\begin{equation}\label{eq:weak.F.condition}
\begin{split}
\|F(u)\|_U \leq & L(\|u\|_U + 1),
\\
\| F'(u) \psi \|_U \leq L \|\psi\|_U, \quad & \|\Lambda^{-\frac{1}{2}} F''(u) (\psi_1, \psi_2) \|_U \leq L \|\psi_1\|_U \|\psi_2\|_U
\end{split}
\end{equation}
for all $u, \psi, \psi_1, \psi_2 \in U$.
Furthermore we assume \eqref{eq:introd.weak.QA} holds for some $\beta \in (0, 1)$.
\end{assumption}
Several comments should be added concerning the above assumption. We first mention that the essential condition \eqref{eq:introd.weak.QA} has been also used in \cite{KLL11,KLL12}, to study weak errors for the linear SWE ($F \equiv 0$) with additive noise. Secondly, it should be pointed out that $F$ in the previous assumption is not required to be twice Fr\'{e}chet differentiable in $U$, which is in general not fulfilled for the Nemytskij operators (see Example \ref{examp:weak.convergence} below). Besides, recall that $\big\| \Lambda^{\frac{\beta -1}{2}} Q^{\frac{1}{2} } \big\|^2_{\mathcal{L}_2(U)} \leq \big\| \Lambda^{ -\frac{1}{2} + \beta} Q  \Lambda^{ -\frac{1}{2} } \big\|_{\mathcal{L}_1(U)}$ (see Theorem 2.1 in \cite{KLL12}). Accordingly, it is easy to see that Assumption \ref{assump:FG} is fulfilled with $\delta = \beta$ under Assumption \ref{ass:weak.error.assump}. In other words, Assumption \ref{ass:weak.error.assump} suffices to guarantee a unique mild solution of \eqref{eq:SEE.additive} and the strong convergence order of $\min(\beta, 1)$ for the underlying scheme \eqref{eq:exp.Euler}. When $\beta \in [1, \infty)$, Theorem \ref{thm:strong.conv.result} shows that the scheme possesses the strong order of 1. This trivially implies weak convergence order of 1, which can not be improved further as shown in the later weak convergence result (Theorem \ref{thm:weak.conv.result}). Based on the above observations, we take condition \eqref{eq:introd.weak.QA} with $\beta \in (0, 1)$ here, instead of $\beta \in (0, \infty)$.

\begin{theorem} \label{thm:weak.conv.result}
Suppose that all conditions in Assumption \ref{ass:weak.error.assump} are satisfied and assume additionally that for $\kappa = \min(\beta, \tfrac{1}{2}- \epsilon)$ with arbitrarily small $\epsilon \in (0, \tfrac{1}{2})$,
\begin{equation}\label{eq:assumption.F'.weak}
\| F'(u) z \|_{-1}
 \leq C \| z \|_{-\kappa} \cdot \big( \|u\|_{\kappa} + 1 \big), \quad u \in \dot{H}^{\kappa},\: z \in U,
\end{equation}
and $X_0 = (u_0, v_0)^T \in H^1$ is deterministic. Then there exists a constant $C_{\epsilon}$ depending on $T, \beta,\epsilon, \Phi, L$ and the initial data, such that
\begin{equation}\label{eq:weak.mainresult}
\big|\mathbb{E} \! \left[ \Phi (X(T))\right] - \mathbb{E} \! \left[ \Phi (X_M) \right] \big| \leq
C_{\epsilon} \tau^{\min( 2\beta, \frac{1}{2} + \beta -
\epsilon, 1) }
\end{equation}
for $\Phi \in \mathcal{C}_b^2 (H; \mathbb{R})$, and where $X(t), t \in [0, T]$ is the mild solution of \eqref{eq:SEE.additive} and $X_m$, $m =0,1,...M$ are produced by the recurrence equation \eqref{eq:Abstr.Exp.Euler}.
\end{theorem}
The proof of Theorem \ref{thm:weak.conv.result} is postponed after some preparations. Here and below, by $\mathcal{C}^k_b(V, V')$ we denote the space of not necessarily bounded mappings from a Banach space $V$ to the other Banach
space $V'$ that have continuous and bounded Fr\'{e}chet derivatives up to order $k$. For the particular case when $V$ is a Hilbert space and $V' = \R$, one can identify the first derivative $D \Psi (X) \in \mathcal{L}(V,\R) $ of a function $\Psi \colon V \rightarrow \R$ with an element in $V$ due to the Riesz representation theorem, i.e.,
\begin{equation}\nonumber
D \Psi (X) \, Z = \left\langle D \Psi (X), Z \right\rangle_V, \quad \forall \quad X, Z \in V,
\end{equation}
and the second derivative $D^2 \Psi (X)$ with a bounded linear operator such that
\begin{equation}\nonumber
D^2 \Psi (X)  (Z_1, Z_2) = \langle D^2 \Psi (X)Z_1, Z_2 \rangle_V, \quad \forall \quad X, Z_1, Z_2 \in V.
\end{equation}
Next, we make some comments concerning the new key condition \eqref{eq:assumption.F'.weak}. Such condition imposed on the derivative operator $F'(u), u \in U$ can be fulfilled in a concrete setting established in Section \ref{sec:examples}. More specifically, we set $U := L^2\big((0,1), \mathbb{R}\big)$ and let $F\colon U \rightarrow U$ be the Nemytskij operator given by $ F (u)(\xi) = f(\xi, u(\xi)), u \in U, \xi \in (0, 1)$ with $f \colon \mathcal{O} \times \mathbb{R} \rightarrow \mathbb{R}$ being a real-valued function. When $f$ has bounded partial derivatives up to order two as required in \eqref{f_condition2}, the condition \eqref{eq:assumption.F'.weak} can be fulfilled (see Example \ref{examp:weak.convergence} below for more details on the verification). Of course, the condition \eqref{eq:assumption.F'.weak} is hard to be satisfied if the associated function $f$ grows super-linearly or is non-Lipschitz.

An immediate consequence of Theorem \ref{thm:weak.conv.result} gives the following result.

\begin{cor}\label{cor:weak.error}
Assume all conditions in Theorem \ref{thm:weak.conv.result} are fulfilled. Then it holds for arbitrarily small $\epsilon \in (0, \tfrac{1}{2})$ that
\begin{equation}\label{eq:weak.error.bound2}
\big|\mathbb{E} \! \left[ \varphi (u(T))\right] - \mathbb{E} \! \left[ \varphi (u_M) \right] \big| \leq C_{\epsilon} \tau^{\min( 2\beta, \frac{1}{2} + \beta -
\epsilon, 1) }
\end{equation}
for $\varphi \in \mathcal{C}_b^2 (U; \mathbb{R})$. Here $u(T)$ and $u_M$ are the first components of $X(T)$ and $X_M$, respectively.
\end{cor}

Comparing the weak error bound \eqref{eq:weak.error.bound2} with the strong one \eqref{eq:strong.conv.additive}, one can find that the rate of weak convergence is, as expected, improved. To see this, let us look at the special interesting case of space-time white noise ($Q = I$). In this case, \eqref{eq:strong.conv.additive} admits a strong convergence order $\beta < \tfrac{1}{2}$, while \eqref{eq:weak.error.bound2} gives a weak convergence order $1 - \epsilon$ for arbitrarily small $\epsilon > 0$. A closely related work on weak convergence of numerical schemes for SWEs can be found in \cite{KLL12}, where only the linear SWE was considered and time discretizations were done by rational approximation to the exponential function. The corresponding weak rates in \cite{KLL12} are $\min( \tfrac{p}{p+1} 2\beta, 1)$, with $p \geq 1$ being method parameters. Particularly, the linear implicit Euler method ($p = 1$) admits weak rate of $\min( \beta, 1)$ and the Crank-Nicolson method ($p = 2$) admits weak rate of $\min( \tfrac{4}{3} \beta, 1)$. Apparently, the exponential scheme attains better weak convergence rates than time discretization schemes in \cite{KLL12}.

To carry out the weak error analysis, we first introduce a function $\mu \colon [0, T] \times H \rightarrow \R$, defined by
\begin{equation}\label{eq:KFunction}
\mu (t,x) = \mathbb{E} \!\left[\Phi (X(t,x))\right], \quad t\in [0,T], \quad x \in H,
\end{equation}
where $\Phi \in \mathcal{C}_b^2(H;\mathbb{R})$ and $X(t,x)$ is the unique mild solution of \eqref{eq:SEE.additive} with the initial value $x \in H$. Owing to \eqref{eq:weak.F.condition}, one can readily infer that
\begin{equation}\label{eq:assump.F'F''}
\begin{split}
\| \BF'(X) Z \|_H  =& \| \Lambda^{-\frac{1}{2}} F'(u) u^z  \|_U \leq \lambda_1^{-\frac{1}{2}} L \| u^z \|_U \leq \lambda_1^{-\frac{1}{2}} L  \|Z\|_H,
\\
\| \BF''(X) (Z_1, Z_2) \|_H = & \|\Lambda^{-\frac{1}{2}} F''(u) ( u^{z_1}, u^{z_2}) \|_U \leq L \|u^{z_1}\|_U \| u^{z_2} \|_U \leq L \|Z_1\|_H  \|Z_2\|_H
\end{split}
\end{equation}
for all $X = (u, v)^T \in H$, $Z = (u^z, v^z)^T \in H$, $Z_1 = (u^{z_1}, v^{z_1})^T \in H$ and $Z_2 = (u^{z_2}, v^{z_2})^T \in H$. This straightforwardly implies $\BF\in \mathcal{C}_b^2(H, H)$. By \cite[Theorem 9.4]{DZ92},  one knows that the solution $X(t,x)$ is twice differentiable with respect to the initial value $x$.
More formally, the process $\zeta^h(t) = \tfrac{\partial X}{\partial x} (t, x) h$ for $t \in [0, T]$, $h \in H$ is the mild solution of the following equation
\begin{equation}\label{eq:Dve}
\begin{split}
\left\{
    \begin{array}{lll} d \zeta^h = (A \zeta^h + \BF'(X(t,x)) \, \zeta^h) \, dt, \quad t \in (0, T], \\
    \zeta^h(0) = h.
    \end{array}\right.
\end{split}
\end{equation}
And $\eta^{h,g}(t) = \tfrac{\partial^2 X}{\partial x^2} (t, x) (h, g) $ for $t \in [0, T]$, $h,g \in H$ is the mild solution of
\begin{equation}\label{eq:D2ve}
\begin{split}
\left\{
    \begin{array}{lll} d \eta^{h,g} = \Big(A \eta^{h,g} + \BF'(X(t,x)) \, \eta^{h,g} + \BF''(X(t,x)) (\zeta^h, \zeta^g) \Big) \, \mbox{d}t, \quad t \in (0, T],
    \\
    \eta^{h,g}(0) = 0.
    \end{array}\right.
\end{split}
\end{equation}
Bearing these facts in mind and differentiating \eqref{eq:KFunction} with respect to $x$ yield
\begin{equation}\label{eq:Dmu}
D \mu (t,x) \, h = \E \! \left[D \Phi (X(t,x)) \, \zeta^h(t)\right], \quad t \in [0, T], \:\: x, h \in H.
\end{equation}
Differentiating \eqref{eq:Dmu} further shows
\begin{equation}\label{eq:D2mu}
D^2\mu(t,x) (h,g) = \E \Big[D^2\Phi (X(t,x)) (\zeta^h(t),\zeta^g(t))\Big] + \E\Big[D\Phi (X(t,x)) \, \eta^{h,g}(t)\Big], \:\: t \in [0, T], \: x, h,g \in H,
\end{equation}
where $\zeta^h(t)$, $\zeta^g(t)$ and $\eta^{h,g}(t)$ for $t \in [0, T], h,g \in H$ are processes defined as above. Then it is well-known (see \cite[Theorem 9.17]{DZ92}, \cite[Theorem 5.4.2]{DZ96}) that $\mu(t,x)$ defined by \eqref{eq:KFunction} is a unique strict solution to the following deterministic PDE:
\begin{equation}\label{eq:kolmogorov}
\begin{split}
\left\{ \!
    \begin{array}{lll} \frac{\partial \mu}{\partial t}(t,x) = \big\langle Ax + \BF(x), D\mu(t,x) \big\rangle_H + \tfrac{1}{2} \text{Tr}\Big[D^2\mu(t,x) B Q^{\frac{1}{2}} \big(BQ^{\frac{1}{2}}\big)^{*} \Big], \:\:\: t \in (0, T], \: x \in\! \mathcal{D}(A),\\
     \mu(0,x) = \Phi(x), \quad x \in H.
    \end{array}\right.
\end{split}
\end{equation}
In order to eliminate the operator $A$, we introduce another process $\nu \colon [0, T] \times H \rightarrow \R$, defined by
\begin{equation}\label{eq:nu}
\nu(t, y) = \mu(t, E(-t)y ), \:\:\: t \in [0, T], \:\: y \in H.
\end{equation}
By virtue of \eqref{eq:kolmogorov}, and also noticing that $\tfrac{\partial}{\partial t} \big(E(-t)y \big) = - A E(-t)y$ for $y \in \mathcal{D}(A)$, one can derive that
\begin{equation}
\begin{split}
\tfrac{\partial \nu}{\partial t}(t,y) =& \tfrac{\partial \mu}{\partial t}(t, E(-t)y) + \big\langle D \mu(t, E(-t)y), - A E(-t) y \big\rangle_H
\\ = &
\big\langle A E(-t)y + \BF(E(-t) y), D\mu(t,E(-t)y) \big\rangle_H
\\ &
+ \tfrac{1}{2} \text{Tr}\Big[D^2\mu(t,E(-t)y) BQ^{\frac{1}{2}} \big(B Q^{\frac{1}{2}}\big)^{*} \Big]
- \big\langle D \mu(t, E(-t)y), \, A E(-t) y \big\rangle_H
\\ = &
\big\langle \BF(E(-t) y),\, D\mu(t,E(-t)y) \big\rangle_H
+ \tfrac{1}{2} \text{Tr}\Big[D^2\mu(t,E(-t)y) B Q^{\frac{1}{2}} \big(B Q^{\frac{1}{2}}\big)^{*} \Big].
\end{split}
\end{equation}
Note that
\begin{equation}\label{eq:Dnu}
D \nu (t,y) z = D \mu(t, E(-t)y) E(-t) z, \quad t \in [0, T],\:\: y, z \in H,
\end{equation}
and
\begin{equation}\label{eq:D2nu}
D^2 \nu (t, y) (z_1, z_2) = D^2\mu(t, E(-t)y )(E(-t) z_1, E(-t) z_2), \quad t \in [0, T],\:\: y, z_1, z_2 \in H.
\end{equation}
The above three identities together with \eqref{ineq4:L1L2} immediately imply
that
\begin{equation}\label{eq:nu.kolmogorov}
\begin{split}
\left\{\!
    \begin{array}{lll} \frac{\partial \nu}{\partial t}(t,y)
    = \tfrac{1}{2} \text{Tr} \Big[D^2\nu(t,y) E(t)B
    Q^{\frac{1}{2}} \big(E(t)B Q^{\frac{1}{2}}\big)^{*}
    \Big]
    + \! \big\langle E(t) \BF( E(-t)y), D\nu(t,y) \big\rangle_H, \:\: t\in (0, T], \: y \in \mathcal{D}(A) \\
     \nu(0,y) = \Phi(y), \quad y \in H.
    \end{array}\right.
\end{split}
\end{equation}
Before proceeding further, we show regularity results on $D \nu (t, y)$ and $D^2 \nu (t, y)$.
\begin{lemma}\label{lemma:Derivatives.X}
Assume Assumption \ref{ass:weak.error.assump} holds. Then $D \nu(t,y) \colon H \rightarrow \R$ given by \eqref{eq:Dnu} and $D^2 \nu(t, y) \colon H \times H \rightarrow \R $ given by \eqref{eq:D2nu} for any $t \in [0, T]$ and $y \in H$ satisfy
\begin{equation}\label{eq:Dnu.D2nu}
\| D \nu (t, y) \|_H \leq C, \quad \mbox{ and } \quad \| D^2 \nu (t, y) \|_{\mathcal{L}(H)} \leq C
\end{equation}
for all $t \in [0, T]$ and $y \in H$.
\end{lemma}
{\it Proof of Lemma \ref{lemma:Derivatives.X}.}
First of all, we shall prove for all $t \in [0, T]$ and $x \in H$ that
\begin{equation}\label{eq:Dmu.D2mu}
\| D \mu (t, x) \|_H \leq C, \quad \mbox{ and } \quad \| D^2 \mu (t, x) \|_{\mathcal{L}(H)} \leq C.
\end{equation}
Then \eqref{eq:Dnu.D2nu} is an easy consequence of \eqref{eq:Dmu.D2mu} thanks to \eqref{eq:Dnu}, \eqref{eq:D2nu} and \eqref{eq:Et.group}. Hence we focus on the proof of \eqref{eq:Dmu.D2mu}. Recall that $\zeta^h (t) = \tfrac{\partial X(t, x)}{\partial x} h $ is the mild solution of \eqref{eq:Dve}, given by
\begin{equation}
\zeta^h(t) = E(t)h + \int_0^t E(t-s)\BF'(X(s,x)) \zeta^h(s) \, \mbox{d}s, \:\: t \in [0, T], \: h \in H.
\end{equation}
Therefore, \eqref{eq:Et.group} and \eqref{eq:assump.F'F''} guarantee that
\begin{equation}\label{eq:zetah}
\|\zeta^h(t)\|_H \leq \|E(t)h \|_H + \lambda_1^{-\frac{1}{2}} L\int_0^t \|\zeta^h(s)\|_H \, \mbox{d}s, \quad t \in [0, T].
\end{equation}
The Gronwall inequality thus shows for all $t \in [0, T]$ that
\begin{equation}\label{eq:zeta^h.bounds}
\| \zeta^h (t) \|_H  \leq C \| h \|_H.
\end{equation}
Since $\Phi \in \mathcal{C}_b^2(H, \R) $ and  \eqref{eq:zeta^h.bounds} holds, one can derive from \eqref{eq:Dmu} that
\begin{equation}\label{eq:Dmu.bound}
| D \mu(t, x) \, h | \leq \E \! \big[ \|D \Phi(X(t, x)) \|_H \cdot \| \zeta^h(t) \|_H \big] \leq C \| h \|_H
\end{equation}
for all $t \in [0, T], h \in H$. Likewise, $\eta^{h,g}(t) = \tfrac{\partial^2 X(t, x)}{\partial x^2}(h, g)$ for $t \in [0, T], h, g \in H$ is given by
\begin{equation}
\eta^{h,g}(t) = \int_0^t E(t-s) \Big(\BF'(X(s,x)) \, \eta^{h,g}(s) + \BF''(X(s,x)) (\zeta^h(s), \zeta^g(s)) \Big) \mbox{d}s.
\end{equation}
With \eqref{eq:assump.F'F''} and \eqref{eq:zeta^h.bounds} at hand, we first obtain
\begin{equation}
\| \eta^{h,g}(t) \|_H \leq C \|h\|_H \cdot\|g \|_H + \lambda_1^{-\frac{1}{2}} L \int_0^t \|\eta^{h,g}(s) \|_H \, \mbox{d}s.
\end{equation}
Once again, we use the Gronwall inequality to get
\begin{equation}
\| \eta^{h,g}(t) \|_H \leq C \|h\|_H \cdot\|g \|_H, \quad t \in [0, T], \:h,g \in H,
\end{equation}
which together with \eqref{eq:D2mu}, \eqref{eq:zeta^h.bounds} implies
\begin{equation}\label{eq:D2mu.bound}
D^2\mu(t,x) (h,g) \leq C \|h\|_H \cdot\|g \|_H, \quad t \in [0, T], \:x, h,g \in H.
\end{equation}
This and \eqref{eq:Dmu.bound} together finish the proof of \eqref{eq:Dmu.D2mu}.
$\square$

As opposed to the parabolic case \cite{AL13,DA10,WG13a}, $D \mu (t, y)$ and $D^2 \mu (t, y)$ defined as above only admit spatial regularity in the space $H$ and this also remains true for $D \nu (t, y)$ and $D^2 \nu (t, y)$. This is due to the lack of the smoothing property of the associated group $E(t)$.
\begin{lemma}\label{lemma:difference.Xtilde}
Suppose that $u_0 \in \dot{H}^1 $, $v_0 \in U$ and Assumption \ref{ass:weak.error.assump} is fulfilled. Then it holds for any $\gamma \in [0, 1)$ and $p \in [2, \infty)$ that
\begin{equation}\label{eq:Holder.regul.nega.space}
\| \tilde{u}(t)- u_m \|_{L^p(\Omega, \dot{H}^{-\gamma})} \leq C \tau^{\min(\beta + \gamma, 1)}
\end{equation}
for $t \in [t_m, t_{m+1}]$, $m = 0, 1, \cdots, M-1$, and where $\tilde{u}(t)$ is the first component of $\tilde{X}(t)$ given by \eqref{eq:X.tilde} and $u_m$ comes from \eqref{eq:exp.Euler}.
\end{lemma}
{\it Proof of Lemma \ref{lemma:difference.Xtilde}.}
Equality \eqref{eq:X.tilde} shows for $t \in [t_m, t_{m+1}]$ that
\begin{equation}\label{eq:Xtilde-X}
\tilde{X}(t) - X_m = \big( E( t - t_m) - \mathcal{I} \big) X_m + E( t - t_m) \Big( \BF (X_m) (t - t_m) + B ( W(t) - W(t_m) ) \Big),
\end{equation}
which implies for $t \in [t_m, t_{m+1}]$ that
\begin{equation}
\begin{split}
\tilde{u}(t)- u_m =& \big( C(t - t_m) - I \big) u_m + \Lambda^{-\frac{1}{2}} S(t- t_m) v_m
\\ &+
\Lambda^{-\frac{1}{2}} S(t- t_m) F(u_m) (t - t_m)
 +
\Lambda^{-\frac{1}{2}} S(t- t_m) ( W(t) - W(t_m)).
\end{split}
\end{equation}
Further, using Lemma \ref{lemma:BDG.inequality}, Proposition \ref{prop:MB.X.Xm}, Corollary \ref{cor:MB.FXmGXm} and  \eqref{eq:SCNts} yields for $t \in [t_m, t_{m+1}]$, $p \in [2, \infty)$ and $\gamma \in [0, 1)$ that
\begin{equation*}
\begin{split}
\big\| \tilde{u}(t)- u_m \big\|^2_{L^p(\Omega, \dot{H}^{-\gamma})}
\leq &
4\big \| \Lambda^{-\frac{\gamma}{2}} \big( C(t - t_m) - I \big) u_m \big\|^2_{L^p(\Omega, U)}
+
4 \big\|\Lambda^{-\frac{1 + \gamma}{2}} S(t- t_m) v_m \big \|^2_{L^p(\Omega, U)}
\\ & +
4(t - t_m)^2 \cdot \big\| \Lambda^{-\frac{1+\gamma}{2}} S(t- t_m) F(u_m) \big\|^2_{L^p(\Omega, U)}
\\ & +
4 C_p^2 (t - t_m) \cdot \big\| \Lambda^{-\frac{1+\gamma}{2}} S(t- t_m) Q^{\frac{1}{2}} \big\|^2_{\mathcal{L}_2(U)}
\\ \leq &
C \tau^{2\min(\beta + \gamma, 1)} + C \tau^{4} \| F(u_m) \|^2_{L^p(\Omega, U)}
\\ & +
C |t - t_m| \cdot  \big\| \Lambda^{-\frac{\beta+\gamma}{2}} S(t- t_m) \big\|^2_{\mathcal{L}(U)} \cdot \big\| \Lambda^{\frac{\beta -1}{2} } Q^{\frac{1}{2}} \big\|^2_{\mathcal{L}_2(U)}
\\ \leq &
C \tau^{2\min(\beta + \gamma, 1)}.
\end{split}
\end{equation*}
The proof of Lemma \ref{lemma:difference.Xtilde} is thus completed. $\square$

Roughly speaking, Lemma \ref{lemma:difference.Xtilde} implies that higher H\"{o}lder regularity in time of the numerical approximation process $\tilde{u}(t)$ can be achieved when measured in the negative Sobolev space $\dot{H}^{-\gamma} = \mathcal{D}(\Lambda^{-\frac{\gamma}{2}})$, $\gamma >0$. This together with the condition \eqref{eq:assumption.F'.weak} is the key ingredient in the following weak error estimates and will be exploited to obtain the expected weak convergence rates (see the estimate of $K_m^{1,2}$ below). Now we are in a position to start the proof of Theorem \ref{thm:weak.conv.result}.

{\it Proof of Theorem \ref{thm:weak.conv.result}.}
First, we define an auxiliary process $\tilde{Y}(t) = E(T -t) \tilde{X}(t)$, $t \in [0, T]$, by
\begin{equation}
\tilde{Y}(t) =  E( T - t_m) X_m + \int_{t_m}^t  E( T - t_m) \BF (X_m) \, \mbox{d} s + \int_{t_m}^t  E( T - t_m) B \, \mbox{d} W(s), \quad t \in [t_m, t_{m+1}].
\end{equation}
The above definition of $\tilde{Y}(t)$ allows for
\begin{equation}
\tilde{Y}(T) = \tilde{X}(T) = X_M \quad \mbox{ and } \quad \tilde{Y}(0) = E(T) X_0.
\end{equation}
This together with \eqref{eq:KFunction} and \eqref{eq:nu} ensures that
\begin{equation}\label{eq:weak.decomposition}
\begin{split}
\E \big[ \Phi( X_M ) \big] - \E \big[ \Phi\big( X(T) \big)\big] = &
\E \big[ \Phi( \tilde{X}(T) ) \big] - \E \big[ \mu(T, X_0) \big]
\\
= &
\E \big[ \Phi( \tilde{Y}(T) ) \big] - \E \big[ \mu(T, X_0) \big]
\\
 = &
\E \big[  \nu(0,\tilde{Y}(T) ) \big] - \E \big[  \nu(T,E(T)X_0 ) \big]
\\
 = &
\E \big[  \nu(0,\tilde{Y}(T) ) \big] - \E \big[  \nu(T,\tilde{Y}(0) ) \big].
\end{split}
\end{equation}
Before further analysis, we define a finite
dimensional subspace $U_N$ of $U$  by
$U_N := \mbox{span} \{e_1, e_2, \cdots, e_N \}$
for $N\in \mathbb{N}$,
and the projection operator $P_N \colon \dot{H}^{\alpha}\rightarrow U_N$ by
$
  P_N\, \phi = \sum_{i=1}^N \langle \phi, e_i \rangle_U e_i, \: \mbox{for} \:\: \phi \in \dot{H}^{\alpha}, \, \alpha \geq -1.
$
Then, we define $\mathcal{P}_N X = (P_N u, P_N v)^T$ for $X = (u, v)^T \in H $. It can be easily verified that $\mathcal{P}_N X \in \mathcal{D}(A)$ for $X \in H$ and that $\lim_{N \rightarrow \infty} \| \mathcal{P}_N X - X \|_H = 0$ for all $X \in H$. At the moment, we apply It\^{o}'s formula \cite[Theorem 4.17]{DZ92} to $\nu(T-t, \mathcal{P}_N \tilde{Y}(t))$
in each interval $[t_m, t_{m+1}]$ and utilize
\eqref{eq:nu.kolmogorov} to obtain
\begin{align}\label{eq:Ito.error.Galerkin}
&\E \big[  \nu(0, \mathcal{P}_N \tilde{Y}(T) ) \big] - \E \big[  \nu(T, \mathcal{P}_N \tilde{Y}(0) ) \big]
\nonumber \\ =&
\sum_{m = 0}^{M-1} \! \left(\E \big[  \nu(T - t_{m+1}, \mathcal{P}_N \tilde{Y}(t_{m+1}) ) \big] - \E \big[  \nu(T - t_{m}, \mathcal{P}_N \tilde{Y}(t_{m}) ) \big] \right)
\nonumber \\  = &
\sum_{m = 0}^{M-1} \bigg( -  \int_{t_m}^{t_{m+1}} \E \!\left[ \tfrac{\partial \nu}{\partial t}(T -t, \mathcal{P}_N \tilde{Y}(t)) \right] \mbox{d}t
\nonumber \\ &  +
\int_{t_m}^{t_{m+1}} \E \!\left[ D \nu (T -t, \mathcal{P}_N \tilde{Y}(t)) \mathcal{P}_N E(T - t_m) \BF( X_m ) \right] \mbox{d}t
\nonumber \\ & +
\frac{1}{2} \E \! \int_{t_m}^{t_{m+1}}  \mbox{Tr} \!\left[ D^2 \nu(T-t, \mathcal{P}_N \tilde{Y}(t) ) \mathcal{P}_N E(T - t_m) B Q^{\frac{1}{2}} \big( \mathcal{P}_N E(T - t_m) B Q^{\frac{1}{2}} \big)^* \right] \mbox{d}t \bigg)
\nonumber \\ = &
\sum_{m = 0}^{M-1} \bigg( \!
\int_{t_m}^{t_{m+1}} \E \!\left[ \Big\langle D \nu (T -t, \mathcal{P}_N \tilde{Y}(t)), \mathcal{P}_N E(T - t_m) \BF( X_m ) - E(T - t) \BF( E(t - T )\mathcal{P}_N \tilde{Y}(t) )\Big\rangle_H \right] \mbox{d}t
\nonumber \\ & +
\frac{1}{2} \E \! \int_{t_m}^{t_{m+1}}  \mbox{Tr} \Big[ D^2 \nu(T-t, \mathcal{P}_N \tilde{Y}(t) ) \Big( \mathcal{P}_N E(T - t_m) B Q^{\frac{1}{2}} \big( \mathcal{P}_N E(T - t_m) B Q^{\frac{1}{2}}\big)^*
\nonumber \\ & \qquad \qquad \qquad -
E(T - t) B Q^{\frac{1}{2}} \big( E(T - t) B Q^{\frac{1}{2}} \big)^* \Big) \Big] \mbox{d}t
\bigg).
\end{align}
Note that $\nu(t, \cdot) \in \mathcal{C}^2_b(H; \R)$ for $t \in [0, T]$ and $\lim_{N \rightarrow \infty}  \mathcal{P}_N X = X $ in $H$.
Therefore, taking limits in \eqref{eq:Ito.error.Galerkin} as $N \rightarrow \infty$ shows
\begin{align}\label{eq:Ito.error}
&\E \big[  \nu(0, \tilde{Y}(T) ) \big] - \E \big[  \nu(T, \tilde{Y}(0) ) \big]
\nonumber \\ =&
\sum_{m = 0}^{M-1} \bigg( \!
\int_{t_m}^{t_{m+1}} \E \!\left[ \Big\langle D \nu (T -t, \tilde{Y}(t)), E(T - t_m) \BF( X_m ) - E(T - t) \BF( E(t - T )\tilde{Y}(t) )\Big\rangle_H \right] \mbox{d}t
\nonumber \\ & +
\frac{1}{2} \E \! \int_{t_m}^{t_{m+1}}  \mbox{Tr} \Big[ D^2 \nu(T-t, \tilde{Y}(t) ) \Big( E(T - t_m) B Q^{\frac{1}{2}} \big( E(T - t_m) B Q^{\frac{1}{2}}\big)^*
\nonumber \\ & \qquad \qquad \qquad -
E(T - t) B Q^{\frac{1}{2}} \big( E(T - t) B Q^{\frac{1}{2}} \big)^* \Big) \Big] \mbox{d}t
\bigg)
\nonumber \\ = &
\sum_{m = 0}^{M-1} \big( K_m^1 + K_m^2 \big).
\end{align}
Now, we further decompose $K_m^1$ and $K_m^2$ as
\begin{equation}\label{eq:Kmq.decomposition}
\begin{split}
K_m^1 =&
\int_{t_m}^{t_{m+1}} \E \!\left[ \Big\langle D \nu (T -t, \tilde{Y}(t)), \, \big( E(T - t_m) - E(T - t)\big) \BF( X_m ) \Big\rangle_H \right] \mbox{d}t
\\ & +
\int_{t_m}^{t_{m+1}} \E \!\left[ \Big\langle D \nu (T -t, \tilde{Y}(t)), \, E(T - t) \big( \BF( X_m ) - \BF( \tilde{X}(t)) \big) \Big\rangle_H \right] \mbox{d}t
\\ = &
K_m^{1,1} + K_m^{1,2},
\end{split}
\end{equation}
and
\begin{align}
K_m^2 =& \frac{1}{2} \E \! \int_{t_m}^{t_{m+1}} \! \mbox{Tr} \Big[ D^2\nu(T-t, \tilde{Y}(t)) \big( E(T - t_m) - E(T-t)\big) B Q^{\frac{1}{2}} \big( E(T - t_m) B Q^{\frac{1}{2}} \big)^* \Big] \mbox{d}t
\nonumber \\ & +
\frac{1}{2} \E \! \int_{t_m}^{t_{m+1}} \! \mbox{Tr} \Big[ D^2\nu(T-t, \tilde{Y}(t))  E(T-t)  B Q^{\frac{1}{2}} \big( ( E(T - t_m) - E(T-t) ) B Q^{\frac{1}{2}} \big)^*  \Big] \mbox{d}t
\nonumber \\ =&
\frac{1}{2} \E \! \int_{t_m}^{t_{m+1}} \! \mbox{Tr} \Big[ D^2\nu(T-t, \tilde{Y}(t)) \big( E(T - t_m) - E(T-t)\big) B Q \big( E(T - t_m) B \big)^* \Big] \mbox{d}t
\nonumber \\ & +
\frac{1}{2} \E \! \int_{t_m}^{t_{m+1}} \! \mbox{Tr} \Big[ D^2\nu(T-t, \tilde{Y}(t))  E(T-t)  B Q \big( ( E(T - t_m) - E(T-t) ) B \big)^*  \Big] \mbox{d}t
\nonumber \\ = &
K_m^{2,1} + K_m^{2,2}.
\end{align}
In \eqref{eq:Kmq.decomposition}, the fact was used that $E(t-T) \tilde{Y}(t) = \tilde{X}(t)$ by the previous definitions. In the sequel, we shall estimate $K_m^{1,1}, K_m^{1,2}, K_m^{2,1}$ and $K_m^{2,2}$ separately. With the aid of
\eqref{eq:MB.FXmFXm}, \eqref{eq:Et.diff} and \eqref{eq:Dnu.D2nu},
we first estimate $K_m^{1,1}$ as follows
\begin{equation}
\begin{split}
\label{eq:Km11}
|K_m^{1,1}| \leq & \int_{t_m}^{t_{m+1}} \E \!\left[ \big\| D \nu (T -t, \tilde{Y}(t)) \big\|_H  \cdot c_1 \tau \big\| \BF( X_m ) \big\|_{H^1} \right] \mbox{d}t
\\ \leq &
C \tau \int_{t_m}^{t_{m+1}}  \E \! \left[ \big\| F(u_m) \big\|_U \right]  \mbox{d} t
\\ \leq &
C \tau^2.
\end{split}
\end{equation}
In order to estimate $K_m^{1,2}$, a common choice is to invoke techniques involved with the Malliavin calculus \cite{AL13,Brehier12,BrehierKopec13,DA10,WG13a}. Here we will provide an alternative way. To be precise, by means of \eqref{eq:Dnu.D2nu}, H\"{o}lder's
inequality, the trigonometric identity and Taylor's formula in Banach space we obtain
\begin{equation}\label{eq:Km12}
\begin{split}
|K_m^{1,2}| \leq & \int_{t_m}^{t_{m+1}} \E \!\left[ \big\| D \nu (T -t, \tilde{Y}(t))\big\|_H \cdot \big\| E(T - t) \big( \BF( X_m ) - \BF( \tilde{X}(t)) \big) \big\|_H \right] \mbox{d}t
\\ \leq &
C  \int_{t_m}^{t_{m+1}} \bigg( \E \! \left[ \big\| \Lambda^{-\frac{1}{2}} S(T - t) \big( F(u_m) - F(\tilde{u}(t)) \big) \big\|_U^2 \right]
\\ & \qquad +
\E \! \left[ \big\| \Lambda^{-\frac{1}{2}} C(T - t) \big( F(u_m) - F(\tilde{u}(t)) \big) \big\|_U^2 \right] \bigg)^{\frac{1}{2}} \mbox{d} t
\\ = &
C  \int_{t_m}^{t_{m+1}} \big\| \Lambda^{-\frac{1}{2}} \big( F(\tilde{u}(t)) - F(u_m)   \big) \big\|_{L^2(\Omega, U)} \, \mbox{d} t
\\ \leq &
C  \int_{t_m}^{t_{m+1}} \! \int_0^1 \big\| \Lambda^{-\frac{1}{2}} F' \big(u_m + r(\tilde{u}(t)- u_m) \big)\cdot \big( \tilde{u}(t)- u_m \big) \big\|_{L^2(\Omega, U)} \, \mbox{d} r \, \mbox{d} t.
\end{split}
\end{equation}
Proposition \ref{prop:MB.X.Xm} ensures that, for $r \in [0, 1]$, $p \in [2, \infty)$, $\kappa = \min(\beta, \tfrac{1}{2} - \epsilon)$ with arbitrarily small $\epsilon >0$,
\begin{equation}
\| u_m + r(\tilde{u}(t)- u_m)\|_{L^p(\Omega, \dot{H}^{\kappa})} \leq \sup_{t \in [0, T]} \| \tilde{u}(t) \|_{L^p(\Omega, \dot{H}^{\kappa})} \leq
\sup_{t \in [0, T]} \| \tilde{X}(t) \|_{L^p(\Omega, H^{\kappa})} \leq C.
\end{equation}
Accordingly, by \eqref{eq:assumption.F'.weak}, \eqref{eq:Holder.regul.nega.space}, H\"{o}lder's inequality and also taking the preceding estimate into account, we derive from \eqref{eq:Km12} that
\begin{align}
\label{eq:Km12.final}
|K_m^{1,2}| \leq & C  \int_{t_m}^{t_{m+1}} \! \int_0^1 \big\| \|  \tilde{u}(t)- u_m \|_{-\kappa} \cdot \big( \| u_m + r(\tilde{u}(t)- u_m) \|_{\kappa} +1 \big) \big\|_{L^2(\Omega, \R)} \, \mbox{d} r \, \mbox{d} t
\nonumber \\ \leq &
C \tau^{1 + \min(\beta + \kappa, 1)} = C \tau^{1 + \min(2\beta, \frac{1}{2} + \beta - \epsilon, 1)}
\end{align}
with $\kappa = \min(\beta, \tfrac{1}{2} - \epsilon)$ for arbitrarily small $\epsilon >0$.
This together with \eqref{eq:Km11} yields
\begin{equation}\label{eq:Km1.final}
|K_m^1| \leq C \tau^{1 + \min(2\beta, \frac{1}{2} + \beta - \epsilon, 1)}.
\end{equation}
Now we turn to the term $K_m^2$. Concerning $K_m^{2,1}$, we use \eqref{eq:Et.group}, \eqref{eq:Dnu.D2nu} and \eqref{eq:introd.weak.QA} to get
\begin{equation}\label{eq:Km21.pre}
\begin{split}
& \Big| \mbox{Tr} \big[ D^2\nu(T-t, \tilde{Y}(t)) \big( E(T - t_m) - E(T-t)\big) B Q \big( E(T - t_m) B \big)^* \big] \Big|
\\ \leq &
\big\| D^2\nu(T-t, \tilde{Y}(t)) \big( E(T - t_m) - E(T-t)\big) B Q \big( E(T - t_m) B \big)^*  \big\|_{\mathcal{L}_1(H)}
\\ \leq &
\big\| D^2\nu(T-t, \tilde{Y}(t)) E(T-t) \big( E(t - t_m) - \mathcal{I} \big) B \Lambda^{\frac{1}{2}-\beta} \big\|_{\mathcal{L}(U, H)}
\\ & \times
\big\| \Lambda^{ -\frac{1}{2} + \beta} Q  \Lambda^{ -\frac{1}{2} } \big\|_{\mathcal{L}_1(U)} \times
\big \| \Lambda^{ \frac{1}{2} } \big( E(T - t_m) B \big)^*  \big\|_{\mathcal{L}(H, U)}
\\ \leq &
C \big\| \big( E(t - t_m) - \mathcal{I} \big) B \Lambda^{\frac{1}{2}-\beta} \big\|_{\mathcal{L}(U, H)} \times \big \| E(T - t_m) B \Lambda^{ \frac{1}{2} }  \big\|_{\mathcal{L}(U, H)},
\end{split}
\end{equation}
where \eqref{eq:SCNts} implies
\begin{equation}
\begin{split}
\big\| ( E(t - t_m) - \mathcal{I} ) B \Lambda^{\frac{1}{2}-\beta} u \big\|^2_H = &  \big\| \Lambda^{-\beta} S(t-t_m) u \big\|_U^2 + \big\| \Lambda^{-\beta} \big(C(t-t_m) - I\big) u \big\|_U^2
\\ \leq &
C (t - t_m)^{2 \min( 2\beta, 1 )} \cdot \|u\|_U,
\end{split}
\end{equation}
and the trigonometric identity gives
\begin{equation}
\begin{split}
\big \| \big( E(T -t_m) B \Lambda^{\frac{1}{2}} \big) u \big\|_H^2 = \big\| S(T-t_m)u \big\|^2_U + \big\| C(T-t_m)u \big\|^2_U = \|u\|_U^2.
\end{split}
\end{equation}
With the aid of these two estimates, we derive from \eqref{eq:Km21.pre} that
\begin{equation}
\Big| \mbox{Tr} \big[ D^2\nu(T-t, \tilde{Y}(t)) \big( E(T - t_m) - E(T-t)\big) B Q \big( E(T - t_m) B \big)^* \big] \Big| \leq C \tau^{\min( 2 \beta, 1 )},
\end{equation}
and hence
\begin{equation}\label{eq:Km21}
|K_m^{2,1}| \leq C \tau^{1 + \min( 2 \beta, 1 )}.
\end{equation}
Following the same arguments as used in the estimate of $K_m^{2,1}$ and noticing that $\| \Lambda^{ -\frac{1}{2} + \beta} Q  \Lambda^{ -\frac{1}{2} } \|_{\mathcal{L}_1(U)} = \| \Lambda^{ -\frac{1}{2} } Q \Lambda^{ -\frac{1}{2} + \beta} \|_{\mathcal{L}_1(U)}$, one can similarly obtain that
\begin{equation}\label{eq:Km22}
|K_m^{2,2}| \leq C \tau^{1 + \min( 2 \beta, 1 )}.
\end{equation}
Gathering \eqref{eq:Km21} and \eqref{eq:Km22} together leads us to
\begin{equation}\label{eq:Km2.final}
|K_m^{2}| \leq C \tau^{1 + \min( 2 \beta, 1 )}.
\end{equation}
Taking \eqref{eq:Km1.final}, \eqref{eq:Km2.final} and \eqref{eq:Ito.error} into account, we conclude form \eqref{eq:weak.decomposition} that
\begin{equation}
\begin{split}
\left| \E \big[ \Phi( X_M ) \big] - \E \big[ X(T) \big] \right| \leq&  \sum_{m = 0}^{M-1} \left( |K_m^1| + |K_m^2| \right)
\\
\leq &
\sum_{m=0}^{M-1} \left(C \tau^{1 + \min(2\beta, \frac{1}{2} + \beta - \epsilon, 1)} + C \tau^{1 + \min(2\beta,1)} \right)
\\ \leq&
C \tau^{\min(2\beta, \frac{1}{2} + \beta - \epsilon, 1)},
\end{split}
\end{equation}
which finishes the proof of Theorem \ref{thm:weak.conv.result}. $\square$

To illustrate the previous assumptions, we give a concrete example as follows.
\begin{example}\label{examp:weak.convergence}
Consider a class of SWEs as discussed in Subsection \ref{subsec:additive.noise}, where we assign here $d = 1$, $\mathcal{O} = (0, 1)$, $g(\xi, u) \equiv 1  $ for all $\xi \in (0, 1), u \in \R$, and $f \colon [0,1] \times \mathbb{R} \rightarrow \mathbb{R}$ is assumed to be a smooth nonlinear function satisfying
\begin{equation}
\label{f_condition2}
|f(\xi, u)| \leq  L(|u|+1),\quad \Big| \tfrac{\partial f}{\partial u}(\xi,u)\Big| \leq L,
\quad \Big| \tfrac{\partial^2 f}{\partial \xi \partial u}(\xi,u)\Big| \leq L, \quad  \mbox{and} \quad  \Big| \tfrac{\partial^2 f}{\partial u^2}(\xi,u)\Big| \leq L
\end{equation}
for all $\xi \in (0,1),\, u\in \mathbb{R}$. Let $U := L^2((0, 1), \R)$ and let
$F\colon U \rightarrow U$ be the Nemytskij operators,
defined by \eqref{eq:Nemytskij}.
Such mappings are, in general,  not Fr\'{e}chet differentiable in $U$, but only G\^{a}teaux
differentiable, with the corresponding derivative operators given by
\begin{align} \label{eq:F.Deriv}
F'(u)(\psi) \,(\xi)= &  \tfrac{\partial f}{\partial u}(\xi,u(\xi)) \cdot \psi(\xi), \quad \xi \in (0,1),
\\
F''(u)(\psi_1, \psi_2) \,(\xi)= &  \tfrac{\partial^2 f}{\partial u^2}(\xi,u(\xi)) \cdot \psi_1(\xi) \cdot \psi_2(\xi), \quad \xi \in (0,1)
\label{eq:F.Deriv2}
\end{align}
for all $u, \psi, \psi_1, \psi_2 \in U$.
Thanks to H\"{o}lder's inequality and a Sobolev inequality: $\dot{H}^{\gamma} = \mathcal{D}(\Lambda^{\frac{\gamma}{2}})$ is continuously embedded into $L^{\infty}((0,1), \R)$ for $\gamma > \tfrac{1}{2}$, we get
\begin{equation}\label{eq:weak.F''}
\begin{split}
\|\Lambda^{-\frac{1}{2}} F''(u) (\psi_1, \psi_2) \|_U
= &
\sup_{\|\psi\|_U \leq 1} \big| \big\langle \Lambda^{-\frac{1}{2}} F''(u) (\psi_1, \psi_2),  \psi \big\rangle_U \big|
\\ = &
\sup_{\|\psi\|_U \leq 1} \big| \big\langle F''(u) (\psi_1, \psi_2),  \Lambda^{-\frac{1}{2}} \psi \big\rangle_U \big|
\\ \leq &
\big\| F''(u) (\psi_1, \psi_2)\|_{L^1((0,1), \R)} \cdot \sup_{\|\psi\|_U \leq 1} \| \Lambda^{-\frac{1}{2}} \psi\|_{L^{\infty}((0,1), \R)}
\\ \leq &
L \|\psi_1\|_{L^2((0,1), \R)} \cdot \| \psi_2 \|_{L^2((0,1), \R)} \cdot C \sup_{\|\psi\|_U \leq 1} \|  \psi\|_U
\\ \leq &
C \|\psi_1\|_U  \| \psi_2 \|_U.
\end{split}
\end{equation}
Now it remains to check \eqref{eq:assumption.F'.weak}.
Since in our example $F$ is a Nemytskij operator, the derivative operator $F'(u), u \in U$ is self-adjoint. Therefore, using the self-adjointness of $\Lambda$ and $F'(u), u \in U$  yields that
\begin{equation}\label{eq:Lambda.F'0}
\begin{split}
\big\| \Lambda^{-\frac{1}{2}} F'(u) z  \big\|_U
= &
\sup_{\| \psi \|_U \leq 1} \big| \big\langle \Lambda^{-\frac{1}{2}} F'(u) z, \psi \big\rangle_U \big|
\\ = &
\sup_{\| \psi \|_U \leq 1} \big| \big\langle z, \big( F'(u)\big)^* \Lambda^{-\frac{1}{2}} \psi \big\rangle_U \big|
\\ = &
\sup_{\| \psi \|_U \leq 1} \big| \big\langle \Lambda^{ -\frac{\kappa}{2} } z, \Lambda^{ \frac{\kappa}{2} }F'(u) \Lambda^{-\frac{1}{2}} \psi \big\rangle_U \big|
\\ \leq &
\big\| \Lambda^{ -\frac{\kappa}{2} } z \big\|_U \cdot
\sup_{\| \psi \|_U \leq 1}  \big\| \Lambda^{ \frac{\kappa}{2} }F'(u) \Lambda^{-\frac{1}{2}} \psi \big\|_U
\end{split}
\end{equation}
with $\kappa = \min(\beta, \tfrac{1}{2} - \epsilon)$ for arbitrarily small $\epsilon >0$. Further,
the setting in Example \ref{examp:weak.convergence} suffices to ensure
\begin{equation}\label{eq:F'.commu.condition}
\| F'(u) \varphi \|_{\gamma} \leq C (\|u\|_{\gamma} + 1) \| \varphi \|_{\sigma}
\end{equation}
for any $u \in \dot{H}^{\gamma}$, $\varphi \in \dot{H}^{\sigma}$, with arbitrarily $\gamma \in (0, \frac{1}{2})$ and $\sigma \in (\frac{1}{2}, 1]$ (see \cite[Lemma 4.4]{WGT13} for more details).
By \eqref{eq:F'.commu.condition} with $\gamma = \kappa = \min(\beta, \tfrac{1}{2} - \epsilon), \sigma =1$ we derive from \eqref{eq:Lambda.F'0} that
\begin{equation}\label{eq:Lambda.F'}
\begin{split}
\big\| \Lambda^{-\frac{1}{2}} F'(u) z  \big\|_U \leq &
\big\| \Lambda^{ -\frac{\kappa}{2} } z \big\|_U \cdot C \sup_{\| \psi \|_U \leq 1} \big( \|u\|_{\kappa} + 1 \big) \| \psi \|_{U}
\\ \leq &
C \| z \|_{-\kappa} \cdot \big( \|u\|_{\kappa} + 1  \big).
\end{split}
\end{equation}
This completes the verification of the key condition \eqref{eq:assumption.F'.weak}. $\square$
\end{example}

\section{Numerical results}
\label{sect:Num.Examp}
In this section, numerical results are included to demonstrate the above assertions. Consider the Sine-Gordon equation subject to noise as follows
\begin{equation}\label{eq:num.examp.SWE}
\left\{
    \begin{array}{lll}
    \frac{\partial^2 u}{\partial t^2} = \frac{\partial^2 u}{\partial \xi^2} - \sin(u) + (\sigma_0 + \sigma_1 u )\dot{W}, \quad  t \in (0, 1], \:\: \xi \in (0, 1),
    \\
     u(0, \xi) = u_0(\xi), \, \frac{\partial u}{\partial t}  (0,\xi) = v_0(\xi), \: \xi \in (0, 1),
     \\
     u(t,0) = u(t,1) = 0,  \: t \in (0, 1],
    \end{array}\right.
\end{equation}
where for simplicity we assume that the covariance
operator $Q$ has the same eigenfunctions as the
Laplacian with Dirichlet boundary condition, i.e.,
$\phi_i = e_i = \sqrt{2} \sin(i \pi \xi), i \in \N$.
Such equation driven by additive noise ($\sigma_1 = 0$)
has been considered as a numerical example in
\cite{CLS13,WGT13}.
Although only semi-discretization in time has been investigated in this article, spatial discretization needs to be done in order to perform the simulations on the computer. To this end we simply spatially discretize \eqref{eq:num.examp.SWE} via a spectral Galerkin method, with $N = 2^{10}$ fixed (see, for example, \cite{WGT13}).  Since the true solutions required in the following are not available, we take numerical solutions produced by the Crank-Nicolson scheme,
using very small stepsize $\tau_{exact} = 2^{-12}$, as reference solutions. Furthermore, the expectations are approximated by averages over 1000 samples in the following numerical tests.

\begin{figure}[htp]
         \centering
         \includegraphics[width=2.8in,height=2.7in]
         {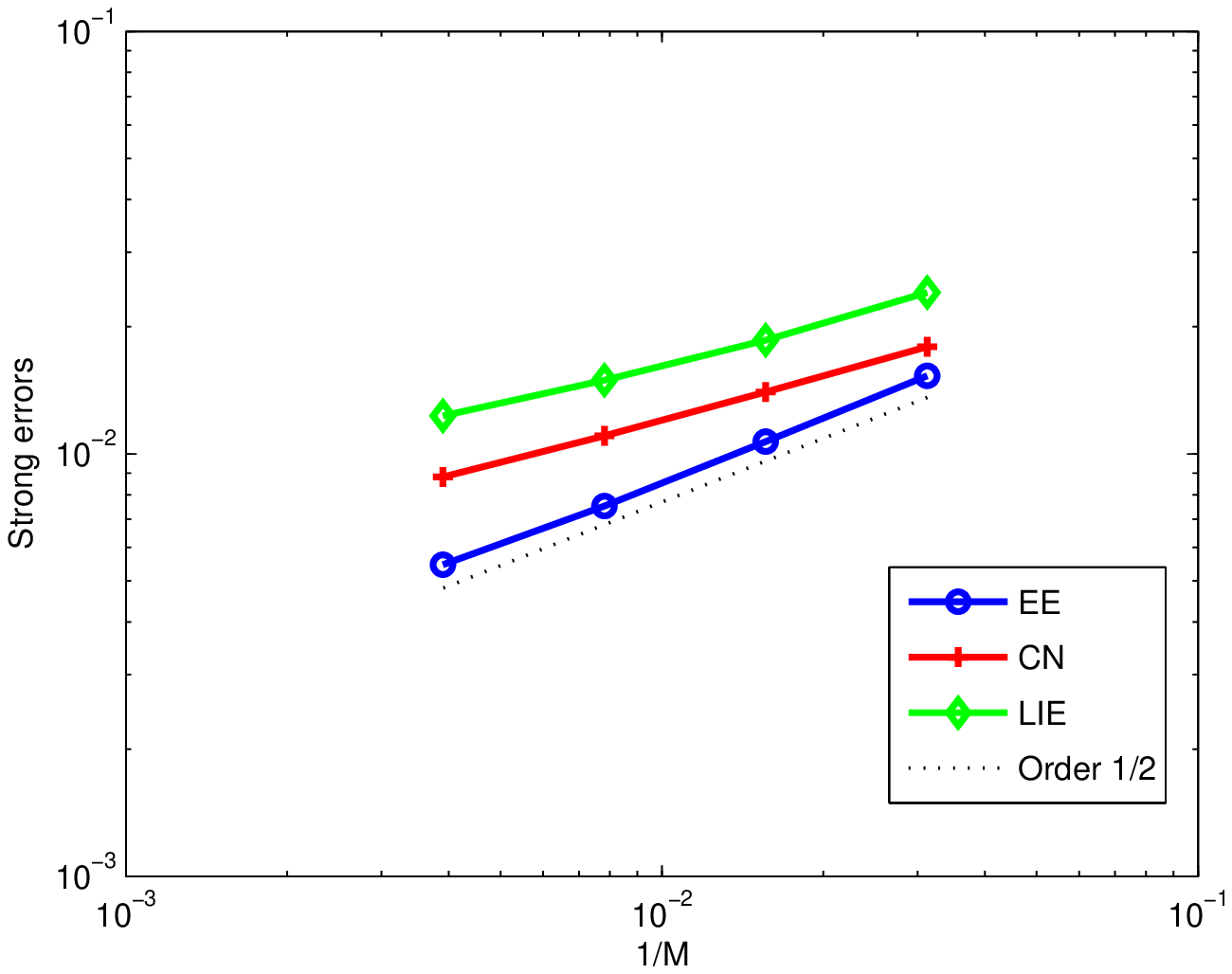}
         \includegraphics[width=2.8in,height=2.7in]
         {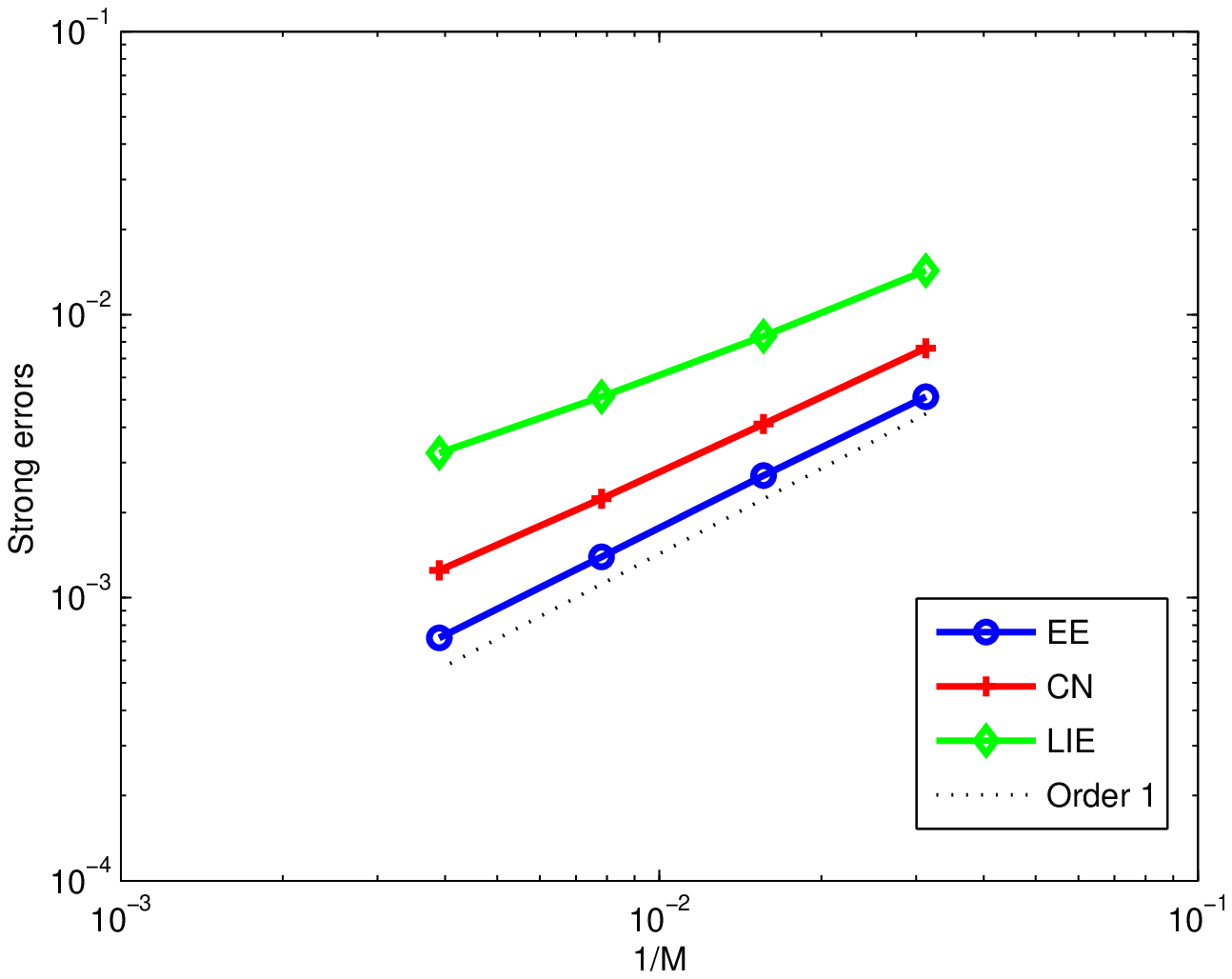}
         \caption{Strong approximation errors for time-stepping schemes applied to \eqref{eq:num.examp.SWE} with space-time white noise (left) and trace class noise (right).}
         \label{fig:strong.error}
\end{figure}

First, let us start with tests on the strong convergence rates and examine the strong approximation error $\left(\E \!\left[ \|u(1) - u_M\|_U^2 \right]\right)^{\frac{1}{2}} = \big( \E \int_0^1 |u(1, \xi) - u_M(\xi)|^2 \, \mbox{d} \xi \big)^{\frac{1}{2}}$, which arises due to the linear implicit Euler (LIE), Crank-Nicolson (CN) and exponential Euler (EE) time discretizations. We choose a set of parameters as $\sigma_0 =0, \sigma_1 = 1$, $u_0(\xi) = 0, v_0(\xi) = \cos(\xi)$. For the case of space-time white noise ($q_i = 1, i \in \N$) and trace-class noise ($q_i = i^{-1.1}, i \in \N$), the left and right plots of Fig.\ref{fig:strong.error}, respectively, depict the strong approximation errors against $\tfrac{1}{M}$ on a log-log scale, with $M = 2^k, k =5,6,7,8$. It turns out that the exponential Euler (EE) scheme \eqref{eq:exp.Euler} performs better than the linear implicit Euler (LIE) and Crank-Nicolson (CN) time discretizations and exhibits the right strong rates, i.e., order $\frac{1}{2}$ for the space-time white noise and order $1$ for the trace-class noise.

\begin{figure}[htp]
         \centering
         \includegraphics[width=3.5in,height=2.5in]
         {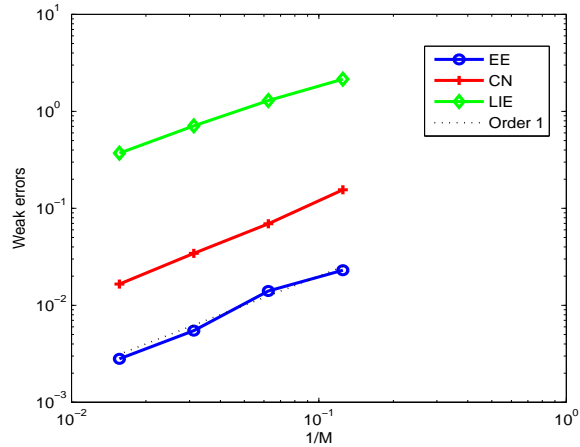}
         \caption{Weak approximation errors for various
 time-stepping schemes
         applied to \eqref{eq:num.examp.SWE} with
space-time white noise.}
         \label{fig:weak.error}
\end{figure}

To illustrate the weak convergence results, we consider the additive space-time white noise and assign $\sigma_0 =1, \sigma_1 = 0$, $u_0(\xi) = \cos(\pi(\xi- \tfrac{1}{2})), v_0(\xi) = 0$. In addition, we choose a particular test function $\varphi(u) = 10\int_0^1 u (\xi) \sin(\pi \xi)\, \mbox{d} \xi$ and consider the weak errors \eqref{eq:weak.error.bound2} when using time-stepping schemes to approximate the quantity  $\E \! \left[ \varphi(u(1)) \right] = 10 \E \big[ \int_0^1 u (1,\xi) \sin(\pi \xi) \, \mbox{d} \xi \big]$. In Table \ref{table:weak.approx}, numerical approximations of $\E \! \left[ \varphi(u(1)) \right] = -4.92169$ by various schemes are presented. Obviously, the scheme \eqref{eq:exp.Euler} always gives more reliable approximations than the other two time-stepping schemes. Even for a large step size $\tau = 2^{-3}$, one can observe a good behavior
of the scheme \eqref{eq:exp.Euler}. Also, weak approximation errors of the three schemes are plotted in Fig. \ref{fig:weak.error}, where one can observe an expected weak rate of $1$ for the scheme \eqref{eq:exp.Euler}.
Despite larger weak approximation errors, the other two schemes (LIE and CN) exhibit
almost weak order one in this
particular example.

\begin{table}[htp]
\begin{center} \footnotesize
\caption{Numerical approximations of $\E \! \left[ \varphi(u(1)) \right] = -4.92169$ by various schemes}
\label{table:weak.approx}
\begin{tabular*}{12cm}{@{\extracolsep{\fill}}cccc}
\hline
$\tau$ & Linear implicit Euler & Crank-Nicolson & Scheme \eqref{eq:exp.Euler} \\
 \hline \\
$2^{-3}$ & -2.77375 & -5.07689 & -4.89865  \\
$2^{-4}$ & -3.62967  &-4.99096 & -4.90763 \\
$2^{-5}$ & -4.21349 & -4.95601 & -4.91620  \\
$2^{-6}$ & -4.55052 & -4.93828 & -4.91888  \\
\hline
\end{tabular*}
\end{center}
\end{table}


\section*{Acknowledgment}
The author thanks the anonymous referee whose insightful comments and valuable suggestions are crucial to the improvement of the manuscript. The author would like to thank Professor Arnulf Jentzen for his financial support and helpful discussions during the author's short visit to ETH Z\"{u}rich in 2013. Thanks also go to Dr. Fengze Jiang for his careful reading the early version of this manuscript. This work was partially supported by National Natural Science Foundations of China under grant numbers 11301550, 11171352, China Postdoctoral Science Foundation under grant numbers 2013M531798, 2014T70779 and Research Foundation of Central South University.

\bibliographystyle{acm}

\end{document}